\documentclass[preprint,review,10pt]{elsarticle}
\biboptions{numbers,sort&compress}

\usepackage[a4paper,left=2.5cm,right=2.5cm,top=3cm,bottom=3cm]{geometry}
\usepackage[parfill]{parskip}  

\usepackage[english]{babel}
\usepackage{lmodern}            
\usepackage{bm}                 

\usepackage{amssymb, amstext, amsthm, mathtools}

\usepackage{graphicx}           
\usepackage{subcaption}         

\usepackage{booktabs, tabularx, float}
\usepackage{enumitem}           

\newcommand{\bx}{\bm{x}}
\newcommand{\bn}{\bm{n}}

\newcommand{\ra}{\mathrm{RA}}
\newcommand{\hs}{\mathrm{HS}}

\newcommand{\pde}{\mathrm{pde}}
\newcommand{\ic}{\mathrm{ic}}

\newcommand{\eps}{\varepsilon}

\newcommand{\spinn}{\mathrm{SPINN}}

\journal{Journal of Computational Physics}

\usepackage[ruled,vlined]{algorithm2e}

\usepackage[bookmarks=true,colorlinks=true,linkcolor=blue,citecolor=red]{hyperref}
\usepackage[nameinlink,capitalize]{cleveref}
\makeatletter
\AtBeginDocument{\def\@citecolor{red}}
\def\ps@pprintTitle{%
  \let\@oddhead\@empty
  \let\@evenhead\@empty
  \def\@oddfoot{}%
  \let\@evenfoot\@oddfoot}
\makeatother

\begin{document}


\begin{frontmatter}

\title{Reliability-Aware Hard--Soft Physics-Informed Neural Networks for Robust Learning of Challenging Partial Differential Equations}

\author[1,2]{Duc Tien Nguyen\fnref{equal}}
\author[3]{Hang Tran\fnref{equal}}
\author[4]{Trinh Minh Tuan}
\author[5]{Nguyen Duc Manh}
\author[6]{Dinh Gia Ninh}


\fntext[equal]{These authors contributed equally to this work.}

\affiliation[1]{%
 organization={College of Engineering and Computer Science, VinUniversity}, 
city={Hanoi},
country={Vietnam}
}

\affiliation[2]{%
organization={Center for AI Research, VinUniversity}, 
city={Hanoi},
country={Vietnam}
}

\affiliation[3]{%
 organization={Department of Computer Science and Engineering, University of North
Texas}, 
city={Texas},
country={USA}
}

\affiliation[4]{%
  organization={Department of Mathematics and Informatics, Hanoi University of Science and Technology}, 
  city={Hanoi},
  country={Vietnam}
}

\affiliation[5]{%
  organization={Department of Mathematics and Informatics, Hanoi University of Science and Technology}, 
  city={Hanoi},
  country={Vietnam}
}

\affiliation[6]{%
 organization={Group of Materials and Structures, School of Mechanical Engineering, Hanoi University of Science and Technology}, 
city={Hanoi},
country={Vietnam}
}

    

\begin{abstract}
Physics-informed neural networks (PINNs) provide a mesh-free framework for solving
partial differential equations, but their training can be limited by loss
imbalance, optimization stiffness, and difficulty in representing localized or
multi-mode solution structures. Hard--soft PINNs (HSPINN) reduce part of this difficulty
by embedding known Dirichlet or periodic constraints directly into the trial
space. However, the resulting fixed admissible representation can still be
poorly conditioned when the free neural component must represent sharp,
localized, or heterogeneous residual structures. This paper develops a
reliability-aware hard--soft PINN (RA-HSPINN) that preserves exact embedded
constraints while introducing a bounded learnable reliability field in the
interior representation. The method is combined with inverse-EMA global loss
balancing and lightweight reliability regularization, while retaining the
standard mean-square residual form. The reliability field is a numerical
modulation variable, not a physical parameter or calibrated probability. We
evaluate RA-HSPINN on nonlinear Burgers equations, periodic convection, a
mixed-boundary Poisson problem, and a mixed first-order Poisson system.
RA-HSPINN reduces the HSPINN relative error by $98.65\%$ for a sharp-gradient
Burgers benchmark, $72.42\%$ for Burgers data with noisy and incompatible
initial conditions, $61.18\%$ for smooth periodic convection, $60.02\%$ for
localized periodic convection, $29.36\%$ for mixed-boundary Poisson, and
$82.17\%$ for a multi-mode mixed first-order Poisson system. An ablation study further shows that the reliability field provides the dominant
recovery mechanism when global loss balancing alone is insufficient. The
results indicate that reliability-aware modulation is most useful for
hard--soft trial spaces that are admissible but difficult to optimize,
especially in localized and multi-mode PDE regimes.
\end{abstract}


\end{frontmatter}


\section{Introduction}

Partial differential equations remain the central mathematical language for computational physics \citep{evans2010partial, leveque2007finite, nguyen2025inverse, nguyen2025modeling, nguyen2025highly}. They describe transport, diffusion, wave propagation, elasticity, fluid flow, electrostatics, and a broad range of coupled multi-physics processes \citep{bartels2016numerical, tonti2014starting}. Classical numerical methods such as finite differences, finite elements, finite volumes, and spectral discretizations are accurate and mature, but they require a discretized computational domain, mesh construction or grid design, stabilization choices, and repeated assembly or time integration \citep{mattiussi1997analysis, teixeira2023finite, moukalled2015finite, dhatt2012finite}. These requirements can become restrictive when the computational domain is irregular, the solution contains localized structures, or the solver must be coupled with sparse measurements, inverse identification, or rapid surrogate evaluation \citep{li2006immersed, karniadakis2021physics}.

Physics-informed neural networks provide an alternative computational paradigm in which the unknown solution is represented by a neural function and the governing equations are imposed through automatic differentiation \citep{raissi2019physics, karniadakis2021physics, tan2024utilizing, nguyen2026st}. Instead of assembling a discrete operator on a mesh, a PINN evaluates pointwise residuals at collocation points and minimizes a composite objective containing the PDE residual and the prescribed boundary or initial conditions \citep{cuomo2022scientific, luo2025physics}. This mesh-free and differentiable formulation is attractive because it can combine data and physics in the same objective and can be adapted to forward, inverse, and data-assimilation problems \citep{cuomo2022scientific}. Nevertheless, the optimization problem induced by the conventional PINN formulation is often difficult. The PDE residual, Dirichlet boundary term, Neumann boundary term, initial-condition term, and optional data terms can differ substantially in magnitude and gradient scale \citep{wang2021understanding,krishnapriyan2021characterizing, wang2022and}. When all constraints are imposed as soft penalties, the optimizer must simultaneously discover a physically valid interior solution and learn the boundary behavior through penalty trade-offs. If the boundary penalty is too small, the boundary is violated; if it is too large, the boundary term can dominate the training dynamics and prevent accurate interior residual minimization \cite{lagaris1998artificial, braga2021self}.

Hard-constrained and hard--soft PINNs were developed to reduce this penalty conflict \citep{lu2021deepxde, sukumar2022exact, sheikholeslami2025physics, zhou2024physics, nguyen2026adaptive}. In a hard--soft PINN (HSPINN), known constraints such as Dirichlet values or periodicity are embedded directly into the neural representation \citep{sheikholeslami2025physics}. For example, a Dirichlet condition can be enforced through a lifting function and a mask that vanishes on the prescribed boundary \citep{lu2021physics}. Periodic conditions can be incorporated by using sinusoidal features \cite{dong2021method, wang2021eigenvector}. Once these constraints are satisfied by construction, the optimizer no longer needs to balance a Dirichlet penalty against the PDE residual. The remaining constraints, such as Neumann fluxes or initial conditions, can still be treated as soft residual terms. This hard--soft viewpoint improves boundary admissibility and often gives a better conditioned optimization problem than a fully soft PINN.

However, exact boundary embedding alone does not solve all training difficulties. A fixed HSPINN still represents the free interior component through a single global neural map. When the solution contains sharp gradients, localized fronts, multi-mode structure, or unreliable data, the hard mask or feature map can improve admissibility while also creating a difficult approximation landscape for the trainable component. This issue is especially visible in nonlinear Burgers problems with high spatial frequencies, advection problems with localized fronts, elliptic problems with multi-mode structures, and training data with noisy or incompatible initial conditions. In these regimes, a hard constraint provides boundary admissibility, but it does not provide a mechanism for spatially or temporally adaptive modulation of the interior representation. This work positions RA-HSPINN as a targeted extension of HSPINN rather than as a replacement for all PINN formulations. The aim is not to make boundary conditions softer, nor to add an ad hoc residual reweighting rule, but to improve the conditioning of the admissible hard--soft trial space when its fixed interior component is too rigid. The method augments the hard--soft representation with a bounded learnable reliability field. For a Dirichlet-constrained problem, the solution is represented as
\begin{equation}
    u_{\theta,\phi}(\bx,t)
    = h_D(\bx,t) + m_D(\bx)\rho_\phi(\bx,t)v_\theta(\bx,t),
    \label{eq:intro_ra}
\end{equation}
where $h_D$ satisfies the prescribed trace, $m_D$ vanishes on the hard boundary, $v_\theta$ is the trainable solution component, and $\rho_\phi\in[\rho_{\min},1]$ is a learnable reliability or modulation field. The method also uses inverse-EMA adaptive global loss weights to balance different residual components during training. The reliability field should not be interpreted as a physical material parameter or a calibrated probability. It is a numerical modulation field that allows the hard--soft neural representation to adapt its local response while preserving the hard constraints embedded by the mask or feature map.

The main contributions are:
\begin{itemize}
    \item We formulate a reliability-aware extension of HSPINNs that preserves exact embedded Dirichlet or periodic constraints while adding a bounded learnable modulation field in the admissible interior representation.
    \item We combine this representation with inverse-EMA global loss balancing and lightweight reliability regularization, yielding a reproducible training procedure that keeps the residual objective in standard mean-square form.
    \item We include comparisons with SPINN, HSPINN, HSPINN with inverse-EMA global weights, and RA-HSPINN to separate the effect of hard embedding, global loss balancing, and reliability-aware modulation.
    \item We provide a product-rule stiffness analysis showing why RA-HSPINN's reliability-aware modulation interacts differently with first-order versus second-order PDE residuals, and demonstrate that casting the HSPINN trial space in a first-order system formulation avoids the additional optimization burden introduced when differentiating the modulated product twice.
\end{itemize}

The remainder of the paper is organized as follows. Section~\ref{sec:method} develops the mathematical formulation, beginning from soft PINNs and hard--soft admissible trial spaces before introducing the reliability-aware representation, adaptive global weighting, reliability regularization, and training algorithm. Section~\ref{sec:setup} gives the experimental setup in sufficient detail to reproduce the numerical results. Section~\ref{sec:results} reports the benchmark results, including a detailed discussion of the mixed first-order Poisson example. Section~\ref{sec:discussion} discusses the numerical behavior, benefits, costs, and limitations of the method, and Section~\ref{sec:conclusion} concludes the paper.

\section{Mathematical Formulation}
\label{sec:method}

This section develops the mathematical structure of RA-HSPINN from the standard soft-penalty PINN formulation. The goal is to separate two ideas that are often mixed in physics-informed learning: admissibility of the trial space and balancing of multiple residual terms. The hard-soft component controls admissibility by embedding known constraints into the neural ansatz. The reliability component enriches the admissible interior representation without relaxing the embedded constraints. The inverse-EMA loss weights then provide a lightweight mechanism for balancing the resulting objective during optimization.

\subsection{General problem setting}
\label{subsec:general_problem}

Let $\Omega\subset\mathbb{R}^{d}$ be a bounded spatial domain with boundary $\partial\Omega$. For time-dependent problems, the space-time domain is denoted by $Q_T=\Omega\times(0,T]$. The unknown field is a scalar function $u(\bx,t)$ governed by a differential equation of the form
\begin{equation}
    \mathcal{N}[u](\bx,t)=f(\bx,t),
    \qquad (\bx,t)\in Q_T,
    \label{eq:generic_pde}
\end{equation}
where $\mathcal{N}[\cdot]$ may be linear or nonlinear and may contain spatial and temporal derivatives. For stationary problems, the time variable is omitted. The boundary is decomposed as
\begin{equation}
    \partial\Omega=\Gamma_D\cup\Gamma_N,
    \qquad \Gamma_D\cap\Gamma_N=\varnothing,
\end{equation}
where $\Gamma_D$ and $\Gamma_N$ denote the Dirichlet and Neumann parts. The auxiliary constraints are written as
\begin{align}
    u(\bx,t)&=g_D(\bx,t), && \bx\in\Gamma_D, \label{eq:general_dirichlet}\\
    \bn\cdot\nabla u(\bx,t)&=g_N(\bx,t), && \bx\in\Gamma_N, \label{eq:general_neumann}\\
    u(\bx,0)&=u_0(\bx), && \bx\in\Omega, \label{eq:general_initial}
\end{align}
with $\bn$ the outward unit normal. Terms that are not relevant to a particular problem, such as an initial condition for a stationary equation, are omitted.

For any neural approximation $\widehat{u}$, the interior residual is
\begin{equation}
    R_{\pde}[\widehat{u}](\bx,t)
    =\mathcal{N}[\widehat{u}](\bx,t)-f(\bx,t).
    \label{eq:generic_residual}
\end{equation}
The Neumann and initial residuals are
\begin{align}
    R_N[\widehat{u}](\bx,t)
    &=\bn\cdot\nabla\widehat{u}(\bx,t)-g_N(\bx,t),
    \qquad \bx\in\Gamma_N, \label{eq:generic_neumann_residual}\\
    R_{\ic}[\widehat{u}](\bx)
    &=\widehat{u}(\bx,0)-u_0(\bx),
    \qquad \bx\in\Omega. \label{eq:generic_ic_residual}
\end{align}
All derivatives are evaluated by automatic differentiation \citep{paszke2017automatic, paszke2019pytorch}. In practice, the expectations below are empirical averages over randomly sampled collocation points.

\subsection{Soft PINN objective and penalty conflict}
\label{subsec:spinn_objective}

A standard soft PINN approximates the solution directly by a neural network $u_\theta(\bx,t)$ and imposes all constraints through penalties \citep{raissi2019physics}. A generic soft objective is
\begin{equation}
\begin{aligned}
    \mathcal{L}_{\spinn}(\theta)
    ={}&\lambda_{\pde}\mathcal{L}_{\pde}(\theta)
    +\lambda_D\mathcal{L}_{D}(\theta)
    +\lambda_N\mathcal{L}_{N}(\theta)
    +\lambda_{\ic}\mathcal{L}_{\ic}(\theta),
\end{aligned}
\label{eq:spinn_loss}
\end{equation}
where
\begin{align}
    \mathcal{L}_{\pde}(\theta)
    &=\mathbb{E}_{Q_T}\left[\left(R_{\pde}[u_\theta]\right)^2\right], \label{eq:spinn_pde_loss}\\
    \mathcal{L}_{D}(\theta)
    &=\mathbb{E}_{\Gamma_D}\left[\left(u_\theta-g_D\right)^2\right], \label{eq:spinn_d_loss}\\
    \mathcal{L}_{N}(\theta)
    &=\mathbb{E}_{\Gamma_N}\left[\left(R_N[u_\theta]\right)^2\right], \label{eq:spinn_n_loss}\\
    \mathcal{L}_{\ic}(\theta)
    &=\mathbb{E}_{\Omega}\left[\left(R_{\ic}[u_\theta]\right)^2\right]. \label{eq:spinn_ic_loss}
\end{align}
The formulation is general, but the optimization problem is a weighted compromise among residual minimization, boundary fitting, and initial-condition fitting. The weights are not known a priori, and the magnitudes and gradients of the partial losses can differ by orders of magnitude. As a result, a small Dirichlet weight can leave boundary errors, whereas a large Dirichlet weight can dominate the interior physics. This penalty conflict motivates the hard-soft construction.

\subsection{Hard-soft admissible trial space}
\label{subsec:hard_soft_space}

The hard-soft formulation embeds exactly known constraints into the neural trial space. For a Dirichlet boundary, the approximation is written as
\begin{equation}
    u_\theta^{\hs}(\bx,t)
    =h_D(\bx,t)+m_D(\bx)v_\theta(\bx,t),
    \label{eq:hspinn_general}
\end{equation}
where $h_D$ is a lifting function and $m_D$ is a mask function satisfying
\begin{equation}
    h_D(\bx,t)=g_D(\bx,t),
    \qquad
    m_D(\bx)=0,
    \qquad \bx\in\Gamma_D.
    \label{eq:hspinn_lift_mask_conditions}
\end{equation}
Therefore,
\begin{equation}
    u_\theta^{\hs}(\bx,t)=g_D(\bx,t),
    \qquad \bx\in\Gamma_D,
    \label{eq:hspinn_exact_bc}
\end{equation}
for every value of $\theta$. The trainable function $v_\theta$ only controls the free component away from the embedded boundary. Thus, the Dirichlet penalty is removed and the hard-soft objective becomes
\begin{equation}
    \mathcal{L}_{\hs}(\theta)
    =\lambda_{\pde}\mathcal{L}_{\pde}^{\hs}(\theta)
    +\lambda_N\mathcal{L}_{N}^{\hs}(\theta)
    +\lambda_{\ic}\mathcal{L}_{\ic}^{\hs}(\theta),
    \label{eq:hs_loss}
\end{equation}
where the partial losses are evaluated using $u_\theta^{\hs}$ instead of $u_\theta$. The omitted Dirichlet term is not approximated or downweighted; it is exactly satisfied by construction.

Periodic structure can also be imposed at the representation level. Instead of adding a penalty that matches opposite spatial boundaries, the network input can be transformed into periodic coordinates. For the periodic transport problems considered later, the feature vector is
\begin{equation}
    \varphi(x,t)=
    \left[
    \sin(x-\beta t),\;\cos(x-\beta t),\;t/t_c
    \right],
    \qquad
    t_c=\frac{1}{\max(\beta,1)}.
    \label{eq:fourier_features}
\end{equation}
Any neural function of $\varphi$ is periodic in the spatial coordinate \(x\) by construction. This mechanism plays the same conceptual role as the mask in Eq.~\eqref{eq:hspinn_general}: known structure is placed in the representation rather than enforced solely by a penalty term. In the mixed first-order periodic Poisson experiment later, the same Fourier idea is used together with a soft boundary/data loss that anchors the manufactured solution on sampled boundary points.

\subsection{Reliability-aware extension of the hard-soft ansatz}
\label{subsec:ra_extension}

\begin{figure}[!t]
\centering
\includegraphics[width=0.98\textwidth]{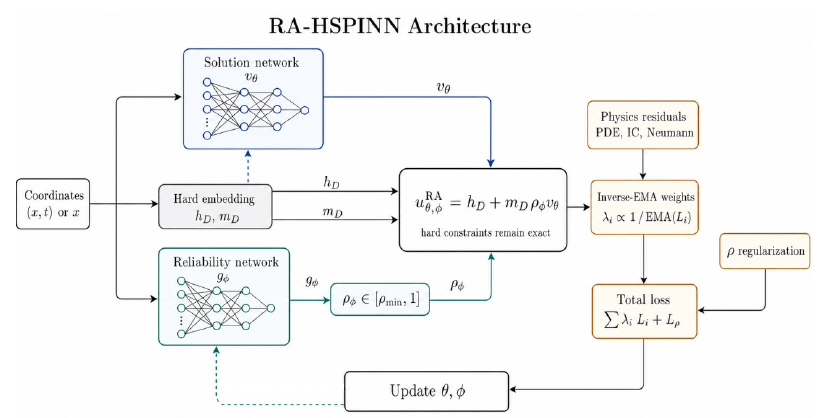}
\caption{Architecture of RA-HSPINN. The solution network \(v_\theta\) is
modulated by a bounded reliability field \(\rho_\phi\), while the hard
embedding \(h_D+m_D(\cdot)\) preserves the prescribed constraints exactly.
Physics residuals are balanced using inverse-EMA adaptive global weights, and
the total objective includes reliability regularization.}
\label{fig:rahspinn_architecture}
\end{figure}

Like HSPINN, RA-HSPINN (see Fig.~\ref{fig:rahspinn_architecture}) restricts the search to an admissible trial space, but the free component $v_\theta$ remains a single global representation. RA-HSPINN augments this free component with a bounded reliability field:
\begin{equation}
    u_{\theta,\phi}^{\ra}(\bx,t)
    =h_D(\bx,t)+m_D(\bx)\rho_\phi(\bx,t)v_\theta(\bx,t).
    \label{eq:ra_general}
\end{equation}
The reliability field is parameterized as
\begin{equation}
    \rho_\phi(\bx,t)
    =\rho_{\min}+(1-\rho_{\min})\sigma(g_\phi(\bx,t)),
    \qquad 0<\rho_{\min}<1,
    \label{eq:rho}
\end{equation}
where $g_\phi$ is a neural network and $\sigma$ is the sigmoid function. The experiments use $\rho_{\min}=0.10$. Because $\rho_\phi$ multiplies only the masked trainable component, the embedded boundary condition remains exact:
\begin{equation}
    u_{\theta,\phi}^{\ra}(\bx,t)=h_D(\bx,t)=g_D(\bx,t),
    \qquad \bx\in\Gamma_D.
    \label{eq:ra_preserves_hard_bc}
\end{equation}
Thus, reliability awareness does not weaken the hard constraint. It only modifies the interior part of the admissible representation.

The field $\rho_\phi$ is a numerical modulation variable, not a physical material parameter and not a calibrated uncertainty probability. It is introduced to give the network a bounded local degree of freedom in regions where the residual structure is difficult. If the hard-soft representation is already sufficient, regularization encourages $\rho_\phi$ to remain close to one and RA-HSPINN approaches the HSPINN limit. If localized residual structures or unreliable constraints make the optimization heterogeneous, the modulation can help the model depart from a fixed global representation while preserving the embedded constraints.

For clarity, the representative forms used in the experiments can be written compactly as follows. A nonzero Dirichlet condition on a one-dimensional interval is enforced by a constant lifting and the mask $x(1-x)$,
\begin{equation}
    u_{\theta,\phi}^{\ra}(x,t)
    =1+x(1-x)\rho_\phi(x,t)v_\theta(x,t).
    \label{eq:burgers_ansatz}
\end{equation}
A two-dimensional homogeneous Dirichlet condition on two coordinate edges is enforced by the mask $xy$,
\begin{equation}
    u_{\theta,\phi}^{\ra}(x,y)
    =xy\rho_\phi(x,y)v_\theta(x,y).
    \label{eq:poisson_ansatz}
\end{equation}
For periodic problems, the hard-soft and RA forms are respectively
\begin{equation}
    u_{\theta}^{\hs}=v_\theta(\varphi),
    \qquad
    u_{\theta,\phi}^{\ra}=\rho_\phi(\varphi)v_\theta(\varphi).
    \label{eq:periodic_ra_ansatz}
\end{equation}
These examples are not separate methods; they are direct instantiations of Eq.~\eqref{eq:ra_general} or the periodic analogue of the same idea.

\subsection{Reliability-modulated residual objective}
\label{subsec:ra_residual_objective}

The reliability-modulated representation changes the hypothesis space, while the residual objective retains the standard mean-square form. Standard HSPINN minimizes
\begin{equation}
    \mathcal{L}_{\pde}^{\hs}
    =\mathbb{E}\left[\left(R_{\pde}[u_\theta^{\hs}]\right)^2\right].
    \label{eq:standard_hs_residual_loss}
\end{equation}
RA-HSPINN evaluates the same residual form on the reliability-modulated approximation:
\begin{equation}
    \mathcal{L}_{\pde}^{\ra}
    =\mathbb{E}\left[
    \left(R_{\pde}[u_{\theta,\phi}^{\ra}](\bx,t)\right)^2
    \right].
    \label{eq:ra_pde_loss}
\end{equation}
Thus, the proposed method does not add an additional pointwise residual reweighting rule. Its adaptivity enters through the admissible representation itself and through the global balancing of the different residual blocks described next. This choice keeps the method closer to the standard PINN objective and avoids introducing a separate spatial reweighting hyperparameter whose effect can depend strongly on the benchmark.

\subsection{Adaptive global loss balancing}
\label{subsec:adaptive_global}

Let the total objective contain $M$ partial losses $\mathcal{L}_1,\ldots,\mathcal{L}_M$, such as the PDE residual, initial-condition residual, Neumann residual, or data residual. At iteration $k$, the scale of each partial loss is tracked by an exponential moving average,
\begin{equation}
    s_i^{(k)}
    =\beta_w s_i^{(k-1)}+(1-\beta_w)\mathcal{L}_i^{(k)},
    \qquad i=1,\ldots,M,
    \label{eq:ema_scale}
\end{equation}
with $\beta_w=0.98$. If $s_i$ is not yet initialized, the current value of $\mathcal{L}_i$ is used. The inverse scale is
\begin{equation}
    \widetilde{\lambda}_i^{(k)}=\frac{1}{s_i^{(k)}+\eps},
    \label{eq:inverse_scale}
\end{equation}
which gives the normalized coefficient
\begin{equation}
    \lambda_i^{(k)}=
    \frac{M\widetilde{\lambda}_i^{(k)}}
    {\sum_{j=1}^{M}\widetilde{\lambda}_j^{(k)}}.
    \label{eq:normalized_lambda}
\end{equation}
The weights are clipped and renormalized as
\begin{equation}
    \lambda_i^{(k)}\leftarrow
    \operatorname{clip}(\lambda_i^{(k)},\lambda_{\min},\lambda_{\max}),
    \qquad
    \lambda_i^{(k)}\leftarrow
    \frac{\lambda_i^{(k)}}{M^{-1}\sum_{j=1}^{M}\lambda_j^{(k)}}.
    \label{eq:lambda_clip_renorm}
\end{equation}
The experiments use $[\lambda_{\min},\lambda_{\max}]=[0.20,5.0]$ and $\eps=10^{-12}$. This inverse-EMA rule is intentionally lightweight. It introduces no extra learnable parameters and provides a reproducible way to prevent persistent domination by one residual component.

\subsection{Reliability regularization and final objective}
\label{subsec:rho_regularization}

The reliability field should not oscillate freely or suppress the solution channel unnecessarily. A regularization term is therefore added:
\begin{equation}
    \mathcal{L}_{\rho}
    =\alpha_{\rho}\mathbb{E}\left[(\rho_\phi-1)^2\right]
    +\alpha_s\mathbb{E}\left[\|\nabla\rho_\phi\|_2^2\right].
    \label{eq:rho_reg}
\end{equation}
The first part biases the method toward the HSPINN limit when modulation is not needed. The second part penalizes rapid variation of the reliability field. In time-dependent one-dimensional problems, the smoothness term is evaluated as $\rho_x^2+0.1\rho_t^2$. In two-dimensional stationary problems, it is evaluated as $\rho_x^2+\rho_y^2$. The reported experiments use $\alpha_{\rho}=10^{-3}$ and $\alpha_s=10^{-5}$.

The final training objective is
\begin{equation}
    \mathcal{L}_{\ra}(\theta,\phi)
    =\sum_{i=1}^{M}\lambda_i\mathcal{L}_i+\mathcal{L}_{\rho}.
    \label{eq:total_ra}
\end{equation}
The set of partial losses depends on the problem class. For time-dependent problems, the losses usually contain the PDE residual and the initial-condition residual. For mixed-boundary stationary problems, they contain the interior residual and the boundary flux residual. For first-order system formulations, they contain the system residual and the boundary or data residual. This notation keeps the method general while allowing each PDE to use its natural residual equations.

Algorithm~\ref{alg:rahspinn} summarizes the complete RA-HSPINN procedure. It is written at the method level rather than for a particular benchmark. The same logic is used for the Burgers, convection, Poisson, and mixed first-order experiments by changing the residual definitions and the applicable constraint losses.

Taken together, RA-HSPINN is built on four components working in concert. The foundation is the HSPINN hard-soft trial space, which satisfies Dirichlet or periodic constraints by construction rather than by penalty. On top of this, a bounded learnable reliability field gives the admissible interior representation the local flexibility it needs without disturbing the embedded constraints. A global inverse-EMA loss balancing rule then keeps the different residual terms from competing destructively during optimization. Finally, a reliability regularization term anchors the modulation field close to its HSPINN limit whenever the problem does not demand strong local adaptation. The experiments in the following section instantiate each of these components across four PDE classes.

\begin{algorithm}[!t]
\small
\setlength{\algomargin}{0.8em}
\caption{Training procedure for RA-HSPINN}
\label{alg:rahspinn}
\DontPrintSemicolon
\SetAlgoLined

\KwIn{
$\mathcal{N}$, $f$, $(h_D,m_D)$ or $\varphi$, collocation sizes,
$\eta$, $K$, $\beta_w$,
$(\lambda_{\min},\lambda_{\max})$.
}
\KwOut{$\theta^\star,\phi^\star$.}

Initialize $v_\theta$, $g_\phi$, and EMA scales $\{s_i\}_{i=1}^{M}$.\;

\For{$k=1,\ldots,K$}{
Sample interior and constraint collocation points.\;

Set
$\rho_\phi=\rho_{\min}+(1-\rho_{\min})\sigma(g_\phi)$.\;

Build
$u_{\theta,\phi}^{\mathrm{RA}}=h_D+m_D\rho_\phi v_\theta$
or
$u_{\theta,\phi}^{\mathrm{RA}}=\rho_\phi v_\theta(\varphi)$.\;

Evaluate by AD:
$R_{\mathrm{pde}}=\mathcal{N}[u_{\theta,\phi}^{\mathrm{RA}}]-f$
and constraint residuals $R_1,\ldots,R_{M-1}$.\;

Form the principal PDE loss
$\mathcal{L}_1=\mathbb{E}[R_{\mathrm{pde}}^2]$.
For first-order systems, add the consistency residuals to this same PDE block.
Set the remaining constraint losses as
$\mathcal{L}_{j+1}=\mathbb{E}[R_j^2]$.\;

Update $s_i$, compute $\lambda_i$, and renormalize them using
Eqs.~\eqref{eq:ema_scale}--\eqref{eq:lambda_clip_renorm}.\;

Compute $\mathcal{L}_{\rho}$ using Eq.~\eqref{eq:rho_reg}.\;

Set
$\mathcal{L}_{\mathrm{RA}}
=\sum_{i=1}^{M}\lambda_i\mathcal{L}_i+\mathcal{L}_{\rho}$.\;

Update
$(\theta,\phi)\leftarrow
(\theta,\phi)-\eta\nabla_{\theta,\phi}\mathcal{L}_{\mathrm{RA}}$.\;
}

\Return{$\theta^\star=\theta,\ \phi^\star=\phi$.}

\end{algorithm}

\section{Experimental Setup}
\label{sec:setup}

The experiments are designed to test RA-HSPINN under several distinct types of numerical difficulty. Burgers examples test nonlinear parabolic dynamics with hard Dirichlet constraints. Convection examples test periodic transport and localized fronts. The mixed-boundary Poisson example tests an elliptic problem with hard Dirichlet and soft Neumann constraints. The mixed first-order Poisson example tests a multi-mode elliptic-like system reformulated so that only first-order derivatives appear in the residual equations.

All runs use PyTorch automatic differentiation. The HSPINN solution network is a fully connected multilayer perceptron with $\tanh$ activations, four hidden layers, and width $64$. The RA-HSPINN solution network has the same architecture as HSPINN, and the reliability network has three hidden layers with width $32$. For periodic convection and mixed first-order Poisson, trigonometric feature maps are used to encode periodic structure. For Burgers and Poisson, raw coordinates are used. The batch sizes and training constants are summarized in \cref{tab:reproducibility}.

\begin{table}[t]
\centering
\caption{Reproducibility settings used in the RA-HSPINN experiments.}
\label{tab:reproducibility}
\resizebox{\textwidth}{!}{%
\begin{tabular}{lcccccc}
\toprule
Problem group & Domain & Steps & Interior points & IC/BC points & Solution network & Reliability network \\
\midrule
Burgers & $(x,t)\in[0,1]\times[0,1]$ & 4000 & 2048 & 512 IC & width 64, depth 4 & width 32, depth 3 \\
Convection & $(x,t)\in[0,2\pi]\times[0,2\pi]$ & 4000 & 2048 & 512 IC & width 64, depth 4 & width 32, depth 3 \\
Poisson mixed BC & $(x,y)\in[0,1]^2$ & 4000 & 2048 & 512 Neumann & width 64, depth 4 & width 32, depth 3 \\
Mixed first-order Poisson & $(x,y)\in[0,2\pi]^2$ & 4000 & 2048 & 512 BC & width 64, depth 4 & width 32, depth 3 \\
\bottomrule
\end{tabular}}
\end{table}

The common first-stage optimizer is Adam with learning rate $10^{-3}$ \citep{kingma2017adammethodstochasticoptimization, zhang2018improved}. The inverse-EMA loss weight parameter is $\beta_w=0.98$, with $\eps=10^{-12}$ and clipping bounds $[0.20,5.0]$. The reliability lower bound is $\rho_{\min}=0.10$. The reliability regularization coefficients are $10^{-3}$ for the mean-to-one term and $10^{-5}$ for the smoothness term. For Burgers with noisy and incompatible initial data, both HSPINN and RA-HSPINN additionally use the same hand-designed initial-condition reliability weight near the boundary corners,
\begin{equation}
    w_{\ic}(x)
    =
    \operatorname{clip}
    \left(
    1-0.75\exp[-(x/\delta)^2]
    -0.75\exp[-((1-x)/\delta)^2],
    0.20,1
    \right),
    \qquad \delta=0.08,
    \label{eq:ic_reliability_weight}
\end{equation}
so that the initial-condition loss is evaluated as
\(\mathcal{L}_{\ic}=\mathbb{E}[w_{\ic}(u_\theta(x,0)-u_0^{\mathrm{train}}(x))^2]\). The width $\delta=0.08$ is chosen to match the spatial extent of the Gaussian-like bumps used to construct the incompatible training data; this weight is a fixed problem-informed prior, not part of the learnable reliability field~$\rho_\phi$, and is applied identically to both HSPINN and RA-HSPINN. Hence the reported comparison isolates the additional effect of the RA reliability field and adaptive global loss balancing. The final relative error is computed on dense evaluation grids using
\begin{equation}
    \varepsilon_{L_2}
    =\frac{\|u_{\mathrm{pred}}-u^*\|_2}{\|u^*\|_2}.
    \label{eq:rel_l2}
\end{equation}
For time-dependent problems, the norm is evaluated over the full space-time grid. For multi-output mixed first-order Poisson, only the solution channel $u$ is used in the reported relative error.

\section{Numerical Results}
\label{sec:results}

\subsection{Burgers equation with sharp gradients}

The first challenging benchmark is a manufactured viscous Burgers equation
\begin{equation}
    u_t+uu_x-\nu u_{xx}=f(x,t),
    \qquad (x,t)\in(0,1)\times(0,1],
    \label{eq:burgers}
\end{equation}
with $\nu=0.01$ and exact solution
\begin{equation}
    u^*(x,t)=1+0.80\sin(4\pi x)\exp(-0.50t).
    \label{eq:burgers_sharp_exact}
\end{equation}
The forcing term is computed analytically from the manufactured solution in
Eq.~\eqref{eq:burgers_sharp_exact}. The boundary values satisfy
$u(0,t)=u(1,t)=1$, and both HSPINN and RA-HSPINN enforce these Dirichlet
conditions exactly through the hard ansatz
\[
    u^{\mathrm{HS}}=1+x(1-x)v_\theta,
    \qquad
    u^{\mathrm{RA}}=1+x(1-x)\rho_\phi v_\theta .
\]
Therefore, the comparison isolates the effect of the reliability-aware interior
modulation under the same admissible boundary-constrained trial space.

This benchmark is difficult for two related reasons. First, the spatial
frequency in Eq.~\eqref{eq:burgers_sharp_exact} is high relative to a simple
smooth Burgers solution. Second, the nonlinear convective term $uu_x$ amplifies
the optimization imbalance between the interior residual and the initial
condition. Moreover, the hard mask $x(1-x)$ can make the free component
$v_\theta$ difficult to approximate near the boundary, because the network must
represent the interior correction through a boundary-vanishing factor.

\begin{figure}[t]
\centering
\includegraphics[width=0.98\textwidth]{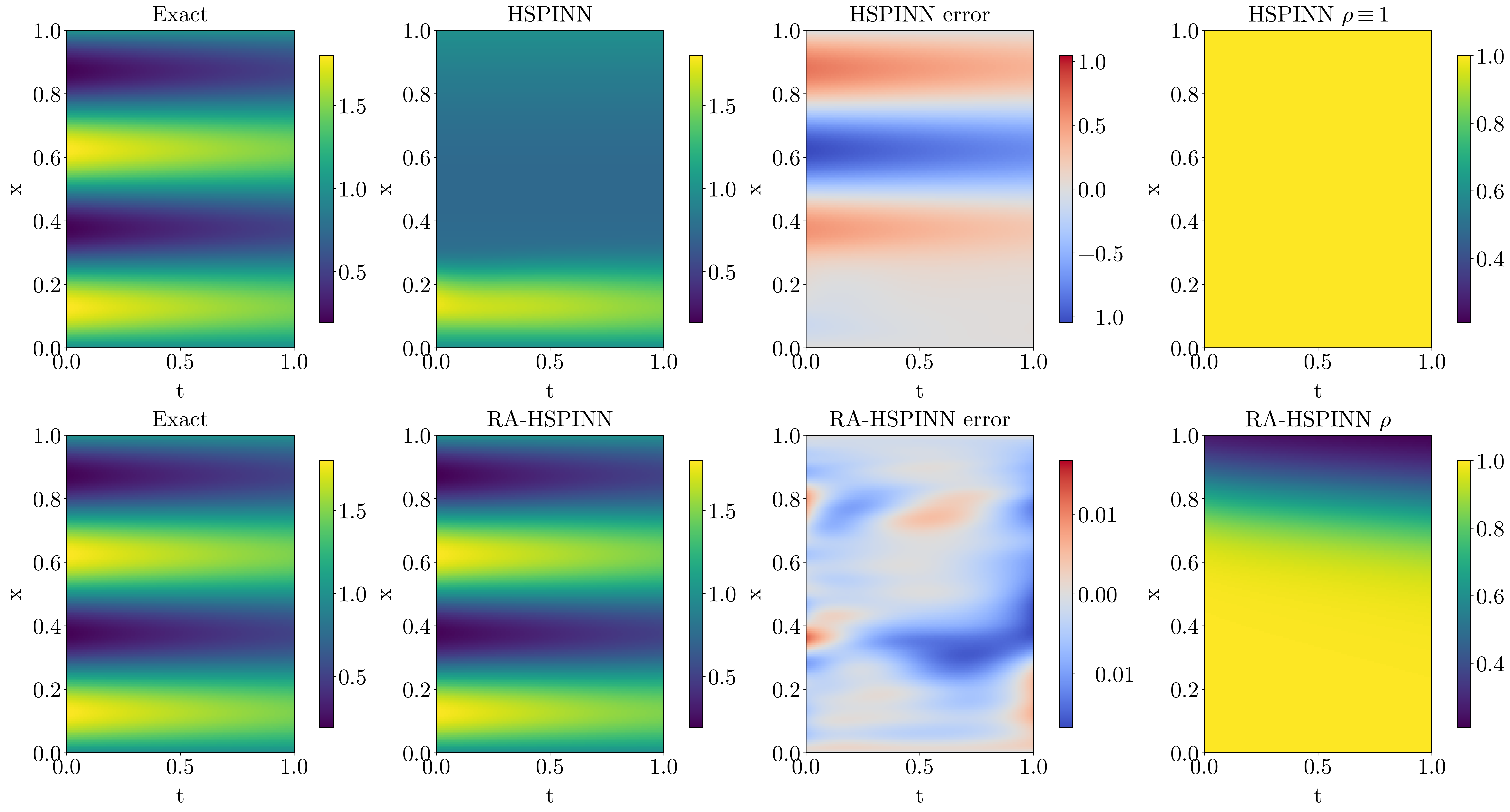}
\caption{Burgers sharp-gradient benchmark. Reference field, HSPINN and
RA-HSPINN predictions, pointwise error fields, and the learned reliability map.
RA-HSPINN suppresses the large space--time error observed in the fixed HSPINN
baseline while preserving the same hard Dirichlet boundary construction.}
\label{fig:burgers_sharp_1}
\end{figure}

The numerical fields in Fig.~\ref{fig:burgers_sharp_1} show that the fixed
HSPINN baseline fails to reproduce the high-frequency Burgers structure and
produces large pointwise errors across the space--time domain. In contrast,
RA-HSPINN closely matches the reference solution and substantially suppresses
the error field. Since both models use the same hard Dirichlet ansatz, the
improvement is not caused by a different treatment of the boundary condition.
Instead, the learned reliability map in Fig.~\ref{fig:burgers_sharp_1}
indicates that the reliability field remains active inside the domain and
locally modulates the free solution component within the same hard-constrained
trial space.

Quantitatively, HSPINN stagnates with a final relative error of
$3.67849\times10^{-1}$, whereas RA-HSPINN reaches $4.965\times10^{-3}$. This
corresponds to a $98.65\%$ reduction in relative error. Both models keep zero
Dirichlet error because the boundary values are embedded exactly by the same
mask. The final mean reliability is approximately $0.838$, confirming that
RA-HSPINN does not simply collapse back to the fixed HSPINN representation.
Rather, it uses a bounded interior modulation to overcome the poor conditioning
of the fixed hard--soft ansatz.

\begin{figure}[t]
\centering
\includegraphics[width=0.65\textwidth]{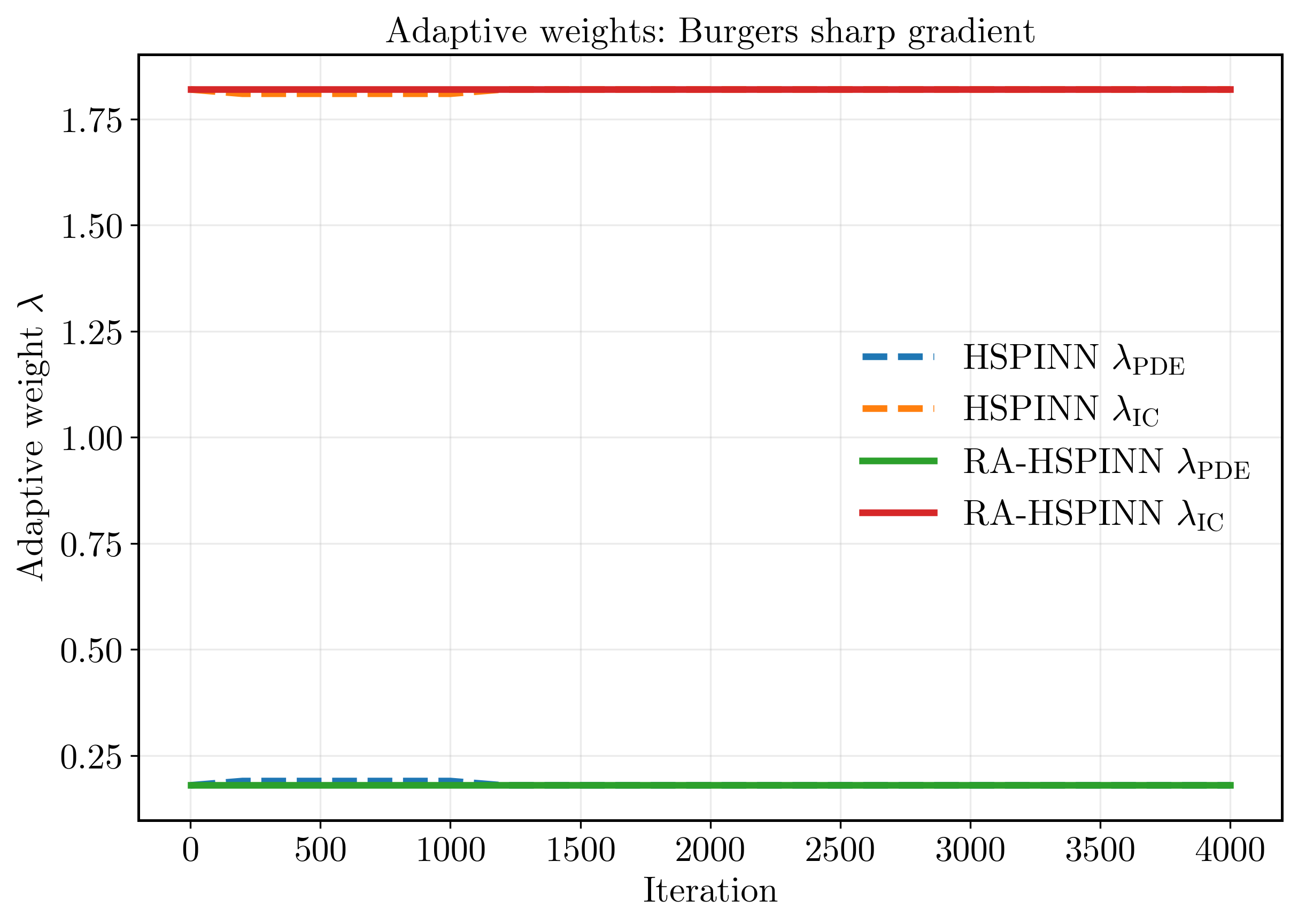}
\caption{Burgers sharp-gradient benchmark. Adaptive global weights for the PDE
and initial-condition losses.}
\label{fig:burgers_sharp_2}
\end{figure}

The training behavior in Fig.~\ref{fig:burgers_sharp_2} further supports this
interpretation. The adaptive global weights balance the PDE and
initial-condition losses during optimization, while the reliability-aware
representation provides the additional local flexibility needed to reduce the
solution error. Thus, this benchmark represents a clear failure mode of the
fixed hard--soft formulation and a strong positive case for the proposed
RA-HSPINN mechanism.

\subsection{Burgers equation with noisy and incompatible initial data}

The second Burgers benchmark tests robustness to unreliable constraints. The
clean reference solution is the smooth manufactured field
\begin{equation}
    u^*(x,t)=1+\sin(\pi x)\exp(-t),
    \label{eq:burgers_clean}
\end{equation}
while the initial condition used during training is deliberately corrupted. The
training initial data are obtained by adding high-frequency sinusoidal noise and
incompatible boundary bumps near $x=0$ and $x=1$. In the implementation, the
high-frequency perturbation is a weighted combination of $\sin(13\pi x)$ and
$\sin(29\pi x)$, and the incompatible component is created using localized
Gaussian-like bumps near both boundary corners. Evaluation is still performed
against the clean exact solution in Eq.~\eqref{eq:burgers_clean}. This design
tests whether the method can avoid overfitting unreliable initial information
while still recovering the clean PDE-governed solution.

Both models use the same hard Dirichlet ansatz, so the boundary values remain
exactly enforced throughout training. In addition, both HSPINN and RA-HSPINN
use the same hand-designed weighting of the unreliable initial-condition data.
Therefore, the comparison isolates the effect of the proposed reliability-aware
interior modulation and adaptive global loss balancing, rather than differences
in boundary treatment or manual data weighting.

\begin{figure}[t]
\centering
\includegraphics[width=0.98\textwidth]{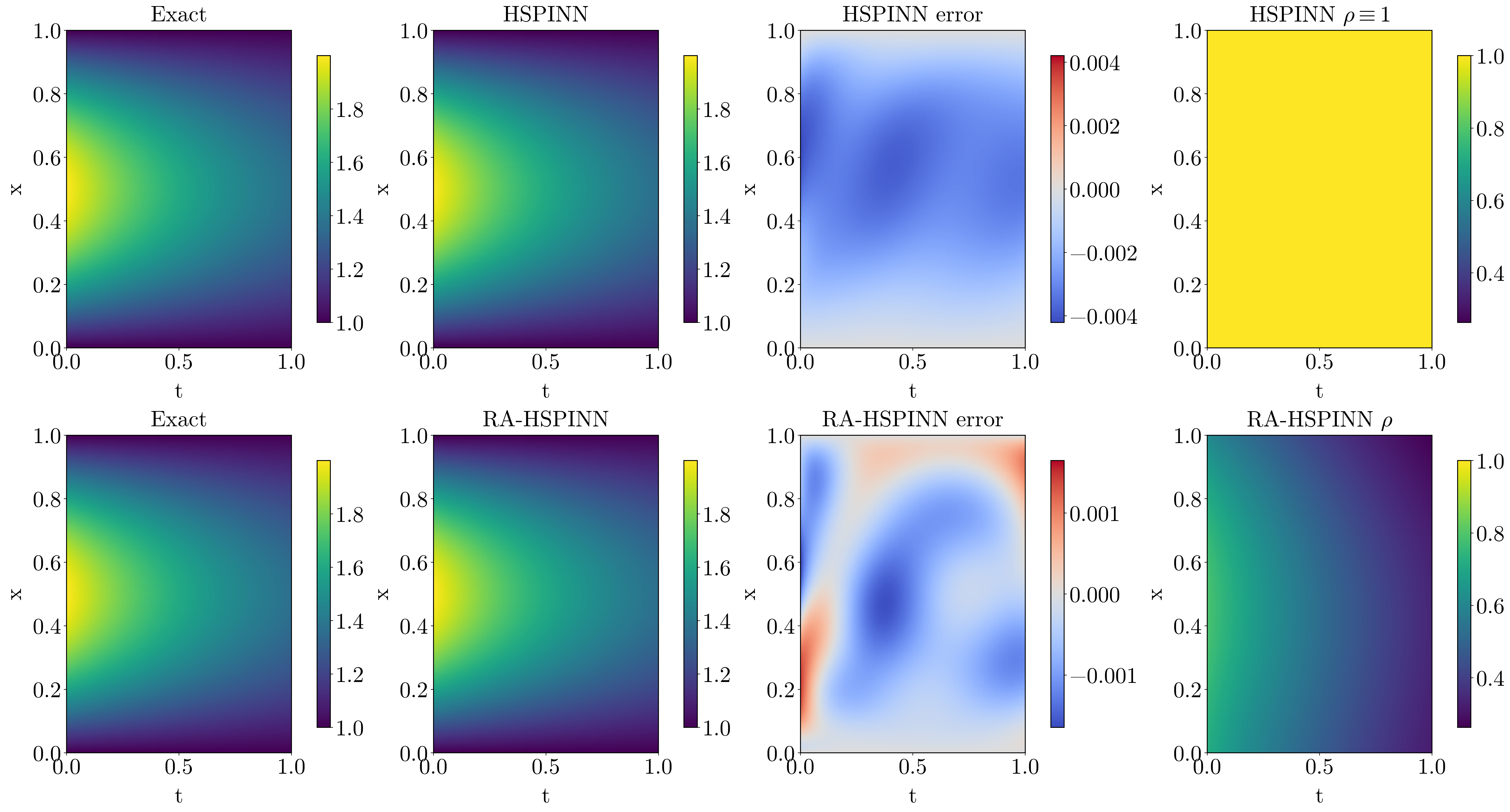}
\caption{Burgers benchmark with noisy and incompatible initial data. Clean
reference field, HSPINN and RA-HSPINN predictions, pointwise error fields, and
the learned reliability map. RA-HSPINN reduces the error while preserving the
same hard Dirichlet boundary construction.}
\label{fig:burgers_noisy_1}
\end{figure}

The numerical fields in Fig.~\ref{fig:burgers_noisy_1} show that both models
recover the clean Burgers solution reasonably well, but RA-HSPINN produces a
visibly smaller error field. The learned reliability map in the same figure
indicates that the modulation is used substantially in this benchmark, which is
consistent with the presence of corrupted and locally incompatible training
constraints. Rather than fitting the noisy initial condition uniformly across
the domain, RA-HSPINN adapts its interior representation to better reconcile the
PDE residual with the unreliable initial data.

Quantitatively, the HSPINN baseline reaches a final relative error of
$1.733\times10^{-3}$, whereas RA-HSPINN reaches $4.78\times10^{-4}$, giving a
$72.42\%$ reduction. The final mean reliability is approximately $0.505$,
showing that the reliability field remains substantially active in this
unreliable-data regime rather than reverting to the fixed HSPINN limit.

\begin{figure}[t]
\centering
\includegraphics[width=0.65\textwidth]{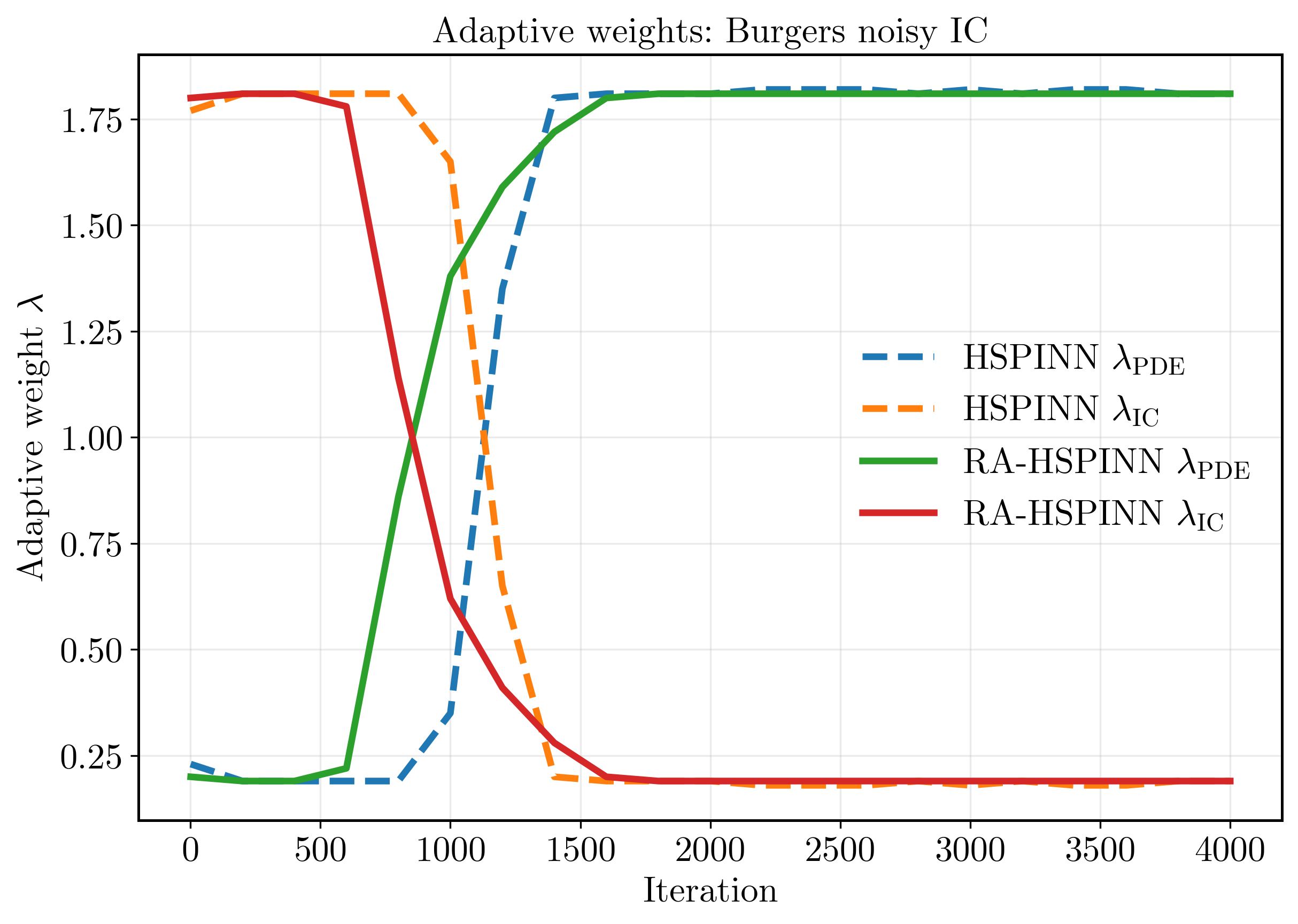}
\caption{Burgers benchmark with noisy and incompatible initial data. Adaptive
global weights for the PDE and initial-condition losses.}
\label{fig:burgers_noisy_2}
\end{figure}

The optimization behavior in Fig.~\ref{fig:burgers_noisy_2} further supports
this interpretation. The adaptive global weights shift strongly toward the PDE
residual by the end of training, with the final RA weights approximately
$\lambda_{\pde}=1.813$ and $\lambda_{\ic}=0.187$. This trend is physically and
numerically consistent with the benchmark design: once the initial condition is
known to contain noise and boundary-incompatible perturbations, the clean PDE
residual should not be dominated by the corrupted data term. Taken together,
Figs.~\ref{fig:burgers_noisy_1}--\ref{fig:burgers_noisy_2} show that RA-HSPINN
improves robustness to unreliable constraints through the joint effect of
reliability-modulated representation and adaptive global loss balancing.

\subsection{Periodic convection with a smooth traveling wave}

The periodic convection equation is
\begin{equation}
    u_t+\beta u_x=0,
    \qquad (x,t)\in[0,2\pi]\times[0,2\pi],
    \label{eq:convection}
\end{equation}
with periodicity in $x$. The smooth traveling-wave solution is
\begin{equation}
    u^*(x,t)=\sin(x-\beta t),
    \label{eq:conv_smooth}
\end{equation}
with $\beta=1$. Periodicity is encoded directly by the feature map in
Eq.~\eqref{eq:fourier_features}; therefore, neither HSPINN nor RA-HSPINN uses a
periodic boundary penalty. The only soft constraints in this benchmark are the
PDE residual and the initial condition.

This problem serves as a clean transport test in which the governing equation
is simple and the exact solution is smooth. Because the periodic structure is
already embedded into the representation through the Fourier feature map, the
baseline HSPINN is expected to perform well. The main question is therefore not
whether RA-HSPINN can repair a severe failure mode, but whether it can still
provide a modest improvement without introducing unnecessary modulation.

\begin{figure}[t]
\centering
\includegraphics[width=0.98\textwidth]{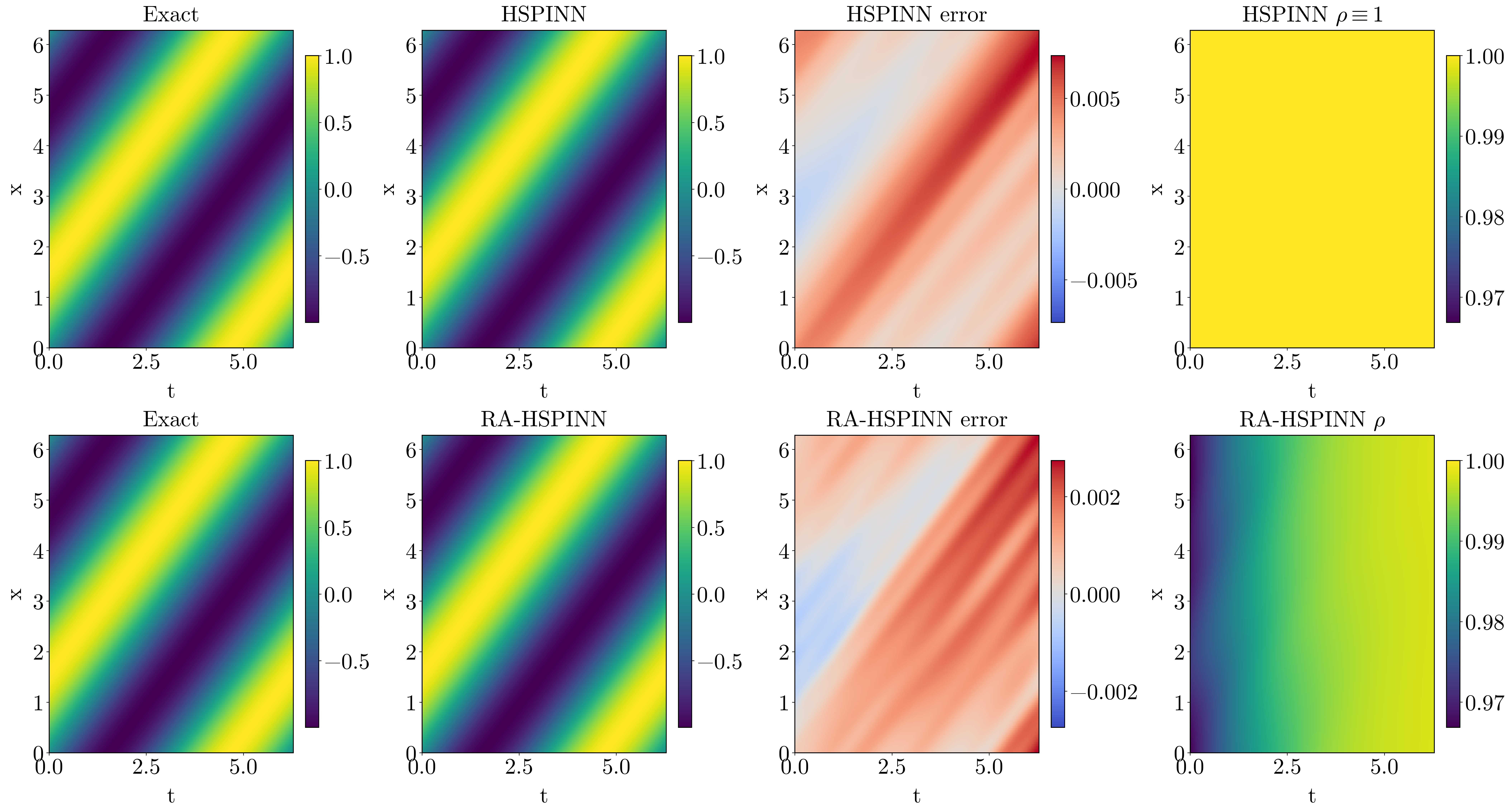}
\caption{Periodic convection benchmark with a smooth traveling wave. Reference
field, HSPINN and RA-HSPINN predictions, pointwise error fields, and the
learned reliability map. Both models capture the periodic transport accurately,
while RA-HSPINN further reduces the error with only mild reliability
modulation.}
\label{fig:fields_convection_beta_1}
\end{figure}

The numerical fields in Fig.~\ref{fig:fields_convection_beta_1} show that both
models reproduce the traveling wave accurately over the full space--time
domain. The HSPINN prediction is already close to the reference field, while
RA-HSPINN further reduces the pointwise error. At the same time, the learned
reliability map in Fig.~\ref{fig:fields_convection_beta_1} remains close to
unity across the domain, indicating that only weak modulation is needed when
the solution is smooth and the periodic feature embedding already provides a
well-conditioned representation.

Quantitatively, HSPINN reaches a final relative error of
$4.432\times10^{-3}$, whereas RA-HSPINN reaches $1.720\times10^{-3}$,
corresponding to a $61.18\%$ reduction. The initial-condition error also
decreases from $4.317\times10^{-3}$ for HSPINN to $8.50\times10^{-4}$ for
RA-HSPINN. However, the final mean reliability is $0.990$, which confirms that
the reliability field is nearly inactive at convergence. In other words, the
improvement does not come from strong local suppression or restructuring of the
solution channel, but from a mild adaptive refinement of an already suitable
hard periodic representation.

\begin{figure}[t]
\centering
\includegraphics[width=0.65\textwidth]{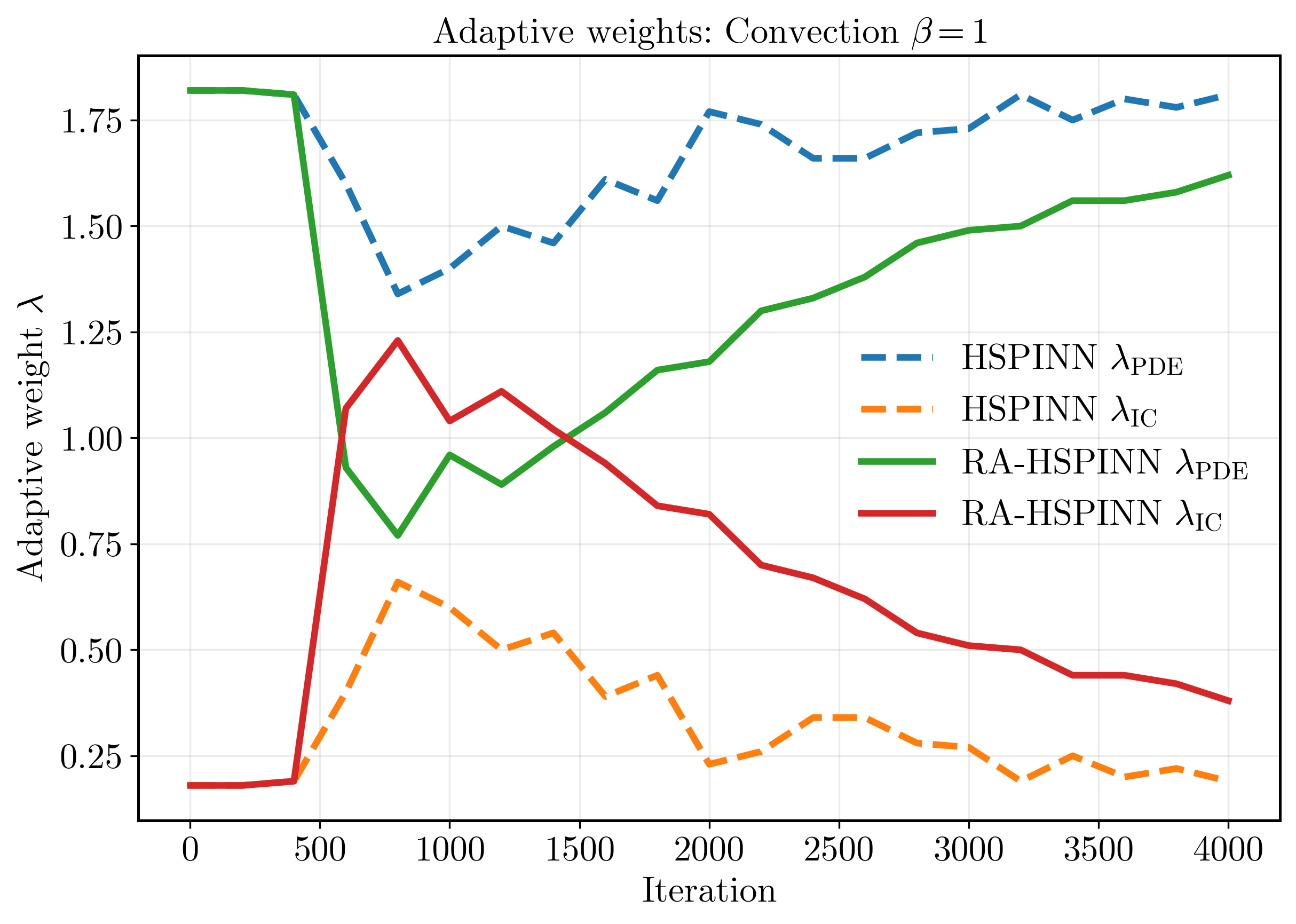}
\caption{Periodic convection benchmark with a smooth traveling wave. Adaptive
global weights for the PDE and initial-condition.}
\label{fig:fields_convection_beta_2}
\end{figure}

The adaptive-weight curves in Fig.~\ref{fig:fields_convection_beta_2} show how
the inverse-EMA rule balances the PDE residual and the initial-condition loss
during training. In both HSPINN and RA-HSPINN, the PDE weight becomes larger
than the initial-condition weight after the early stage, indicating that the
optimization increasingly emphasizes the transport residual once the initial
wave profile has been fitted. For RA-HSPINN, the initial-condition weight
temporarily increases during the early iterations and then decreases as the PDE
weight becomes dominant. This behavior is consistent with the smooth
traveling-wave setting: the initial condition provides an informative anchor,
while the PDE residual controls the propagation of the periodic wave over the
space--time domain. Together with the field and reliability-map results in
Fig.~\ref{fig:fields_convection_beta_1}, these curves show that RA-HSPINN
improves the smooth convection solution through mild reliability-aware
refinement within an already suitable periodic representation, rather than
through strong reliability modulation.

\subsection{Periodic convection with a steep localized front}

The fourth benchmark keeps the same convection equation but replaces the smooth
sine wave with a localized front-like profile,
\begin{equation}
    u^*(x,t)=\exp\left(\kappa[\cos(x-\beta t)-1]\right),
    \qquad \beta=1,
    \qquad \kappa=15.
    \label{eq:steep_front}
\end{equation}
This solution remains an exact solution of the linear advection equation, but
it contains a concentrated peak and a large dynamic range. It is therefore a
useful benchmark because the governing PDE is simple while the solution
representation is locally demanding. A method that relies too strongly on a
global average residual may capture the broad background transport but still
smooth out the localized front.

As in the smooth traveling-wave case, periodicity is embedded directly through
the Fourier feature map, so neither HSPINN nor RA-HSPINN uses a periodic
boundary penalty. The only soft constraints are the PDE residual and the
initial condition. The key question in this benchmark is whether the
reliability-aware formulation can improve the representation of a localized
transport structure even when the underlying PDE is linear.

\begin{figure}[t]
\centering
\includegraphics[width=0.98\textwidth]{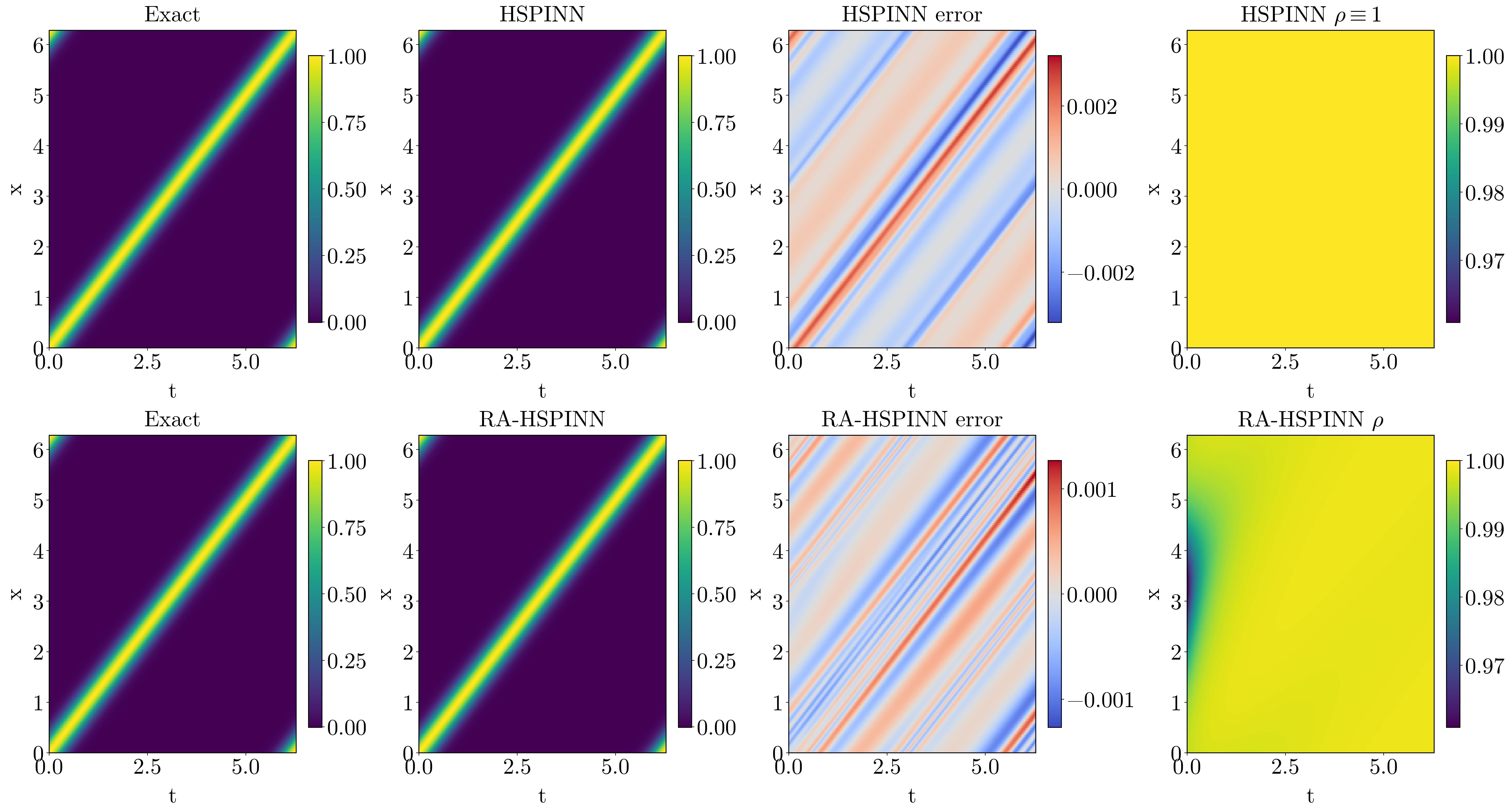}
\caption{Periodic convection benchmark with a steep localized front. Reference
field, HSPINN and RA-HSPINN predictions, pointwise error fields, and the
learned reliability map. RA-HSPINN improves the reconstruction of the
localized peak while the periodic structure is imposed through the feature map.}
\label{fig:fields_convection_steep_front}
\end{figure}

The numerical fields in Fig.~\ref{fig:fields_convection_steep_front} show that
both models reproduce the moving front reasonably well, but RA-HSPINN yields a
clearer reconstruction of the localized peak and a smaller pointwise error
field. The learned reliability map in the same figure remains close to one over
most of the domain, indicating that the periodic feature embedding already
captures the dominant transport structure. Nevertheless, the reliability field
still provides mild local modulation, which is sufficient to improve the
representation of the steep front.

Quantitatively, HSPINN reaches a final relative error of
$3.309\times10^{-3}$, whereas RA-HSPINN reaches $1.323\times10^{-3}$, giving a
$60.02\%$ reduction. RA-HSPINN also reduces the maximum initial-condition error
from $2.27\times10^{-3}$ to $6.71\times10^{-4}$. The final reliability
statistics are mean $0.998$, minimum $0.965$, and maximum $0.999$, confirming
that only mild modulation is needed in this clean advection case. Thus, the
benefit of RA-HSPINN here is not a large suppression effect, but an additional
adaptive flexibility around the localized front within an already suitable hard
periodic representation.

\begin{figure}[t]
\centering
\includegraphics[width=0.65\textwidth]{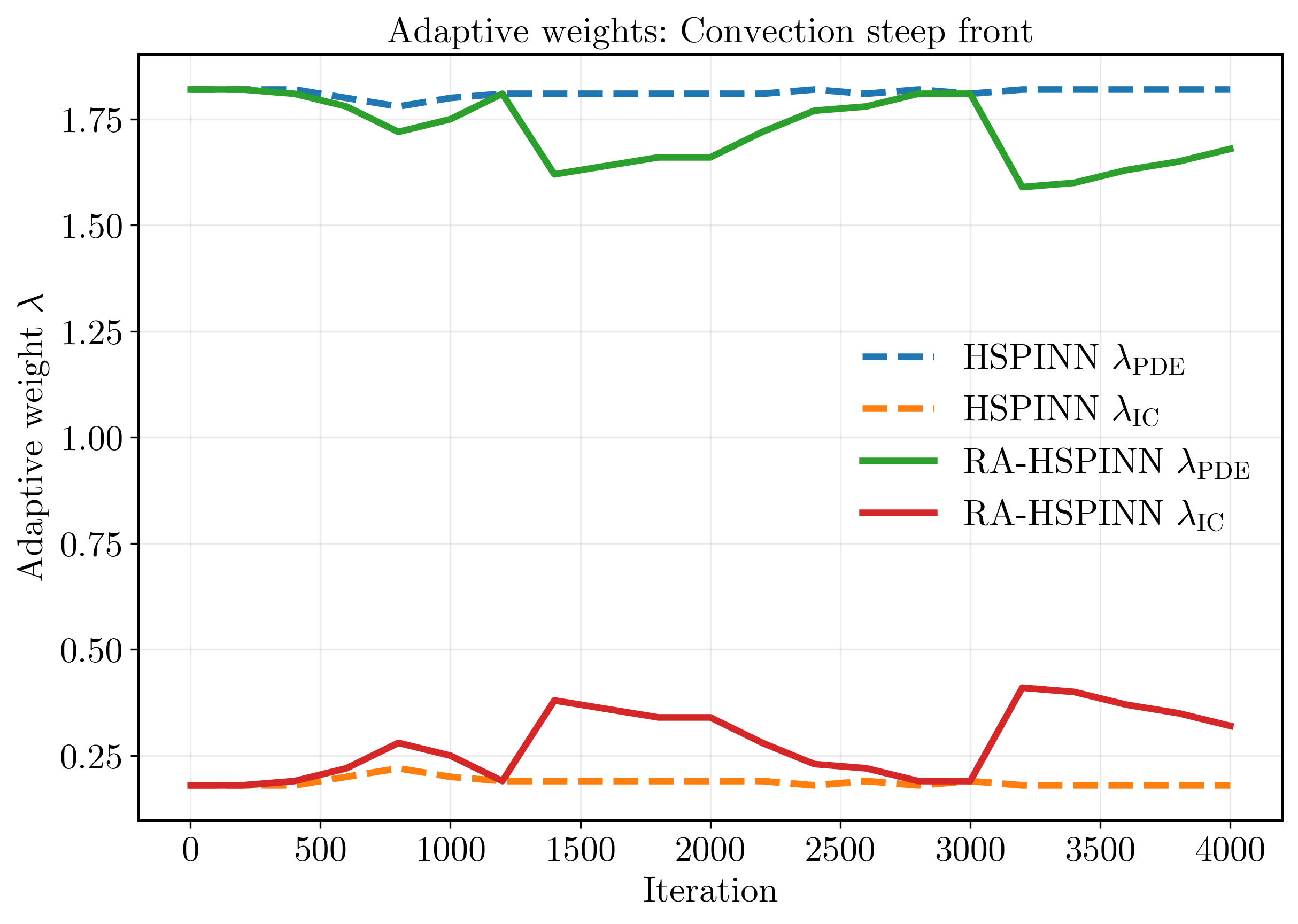}
\caption{Periodic convection benchmark with a steep localized front. Adaptive
global weights for the PDE and initial-condition losses.}
\label{fig:fields_convection_steep_front_2}
\end{figure}

The adaptive-weight curves in Fig.~\ref{fig:fields_convection_steep_front_2}
show how the inverse-EMA rule balances the PDE residual and the
initial-condition loss for the localized-front problem. For both HSPINN and
RA-HSPINN, the PDE weight becomes larger than the initial-condition weight after
the early training stage, indicating that the optimization increasingly
emphasizes the transport residual once the initial localized profile has been
anchored. In RA-HSPINN, the initial-condition weight briefly increases during
the early iterations and then decreases as the PDE weight rises, while the
reliability field provides mild local modulation around the transported peak.
This behavior is consistent with the structure of the benchmark: the initial
condition defines the localized front, whereas the PDE residual controls its
space--time propagation without smoothing out the peak. Together with the field
and reliability-map results in
Fig.~\ref{fig:fields_convection_steep_front}, these curves show that RA-HSPINN
improves localized transport accuracy through the combination of periodic
feature embedding, adaptive residual balancing, and mild reliability-aware
refinement.

\subsection{Mixed-boundary Poisson equation}

The mixed-boundary Poisson benchmark is posed on the unit square:
\begin{equation}
    -\Delta u=f(x,y),
    \qquad (x,y)\in(0,1)^2,
    \label{eq:poisson}
\end{equation}
with exact solution
\begin{equation}
    u^*(x,y)=\sin(\pi x)\sin(\pi y).
    \label{eq:poisson_exact}
\end{equation}
Homogeneous Dirichlet conditions are imposed at $x=0$ and $y=0$, and Neumann
conditions are imposed at $x=1$ and $y=1$. The hard mask $xy$ enforces the two
Dirichlet edges exactly, while the soft objective contains the Poisson residual,
the Neumann residual, and a small solution regularization term.

This benchmark tests the proposed method in a smoother elliptic setting than
the Burgers examples. Because the exact solution is regular and the Dirichlet
constraints are already embedded exactly through the hard mask, the HSPINN
baseline is expected to perform reasonably well. The main question is therefore
whether reliability-aware modulation can still provide a useful refinement when
the hard--soft representation is already well aligned with the target solution.

\begin{figure}[t]
\centering
\includegraphics[width=0.98\textwidth]{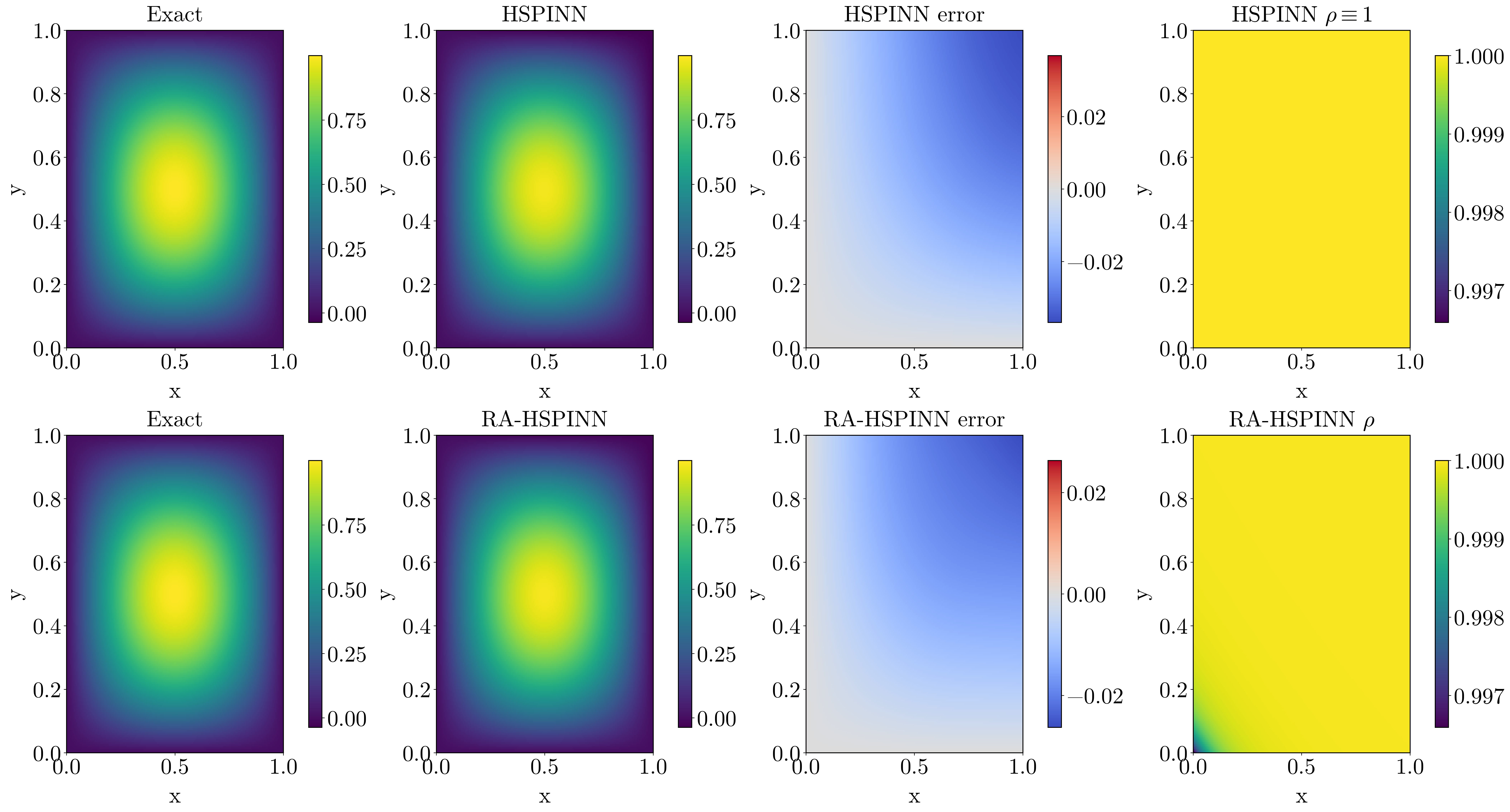}
\caption{Mixed-boundary elliptic Poisson benchmark. Reference field, HSPINN and
RA-HSPINN predictions, pointwise error fields, and the learned reliability map.
Both models satisfy the Dirichlet edges exactly through the hard mask $xy$,
while RA-HSPINN yields a modest reduction in the interior error.}
\label{fig:fields_poisson_base}
\end{figure}

The numerical fields in Fig.~\ref{fig:fields_poisson_base} show that both
models recover the main structure of the Poisson solution, while RA-HSPINN
produces a modestly smaller error field. Since both models use the same hard
mask, both have zero Dirichlet error by construction. The learned reliability
map in Fig.~\ref{fig:fields_poisson_base} indicates that the modulation remains
mild, which is consistent with the absence of the severe conditioning
difficulty observed in the sharp-gradient Burgers problem.

The mixed-boundary Poisson experiment uses the Adam stage followed by the
LBFGS refinement described in \cref{sec:setup}. Quantitatively, HSPINN obtains
a final relative error of $3.5212\times10^{-2}$, whereas RA-HSPINN obtains
$2.4872\times10^{-2}$, corresponding to a $29.36\%$ reduction. This improvement
is meaningful but more modest than in the sharp-gradient Burgers and mixed
first-order Poisson examples. The most natural interpretation is that the exact
solution is smooth and the fixed HSPINN representation is already sufficiently
expressive for most of the approximation task.

\begin{figure}[t]
\centering
\includegraphics[width=0.65\textwidth]{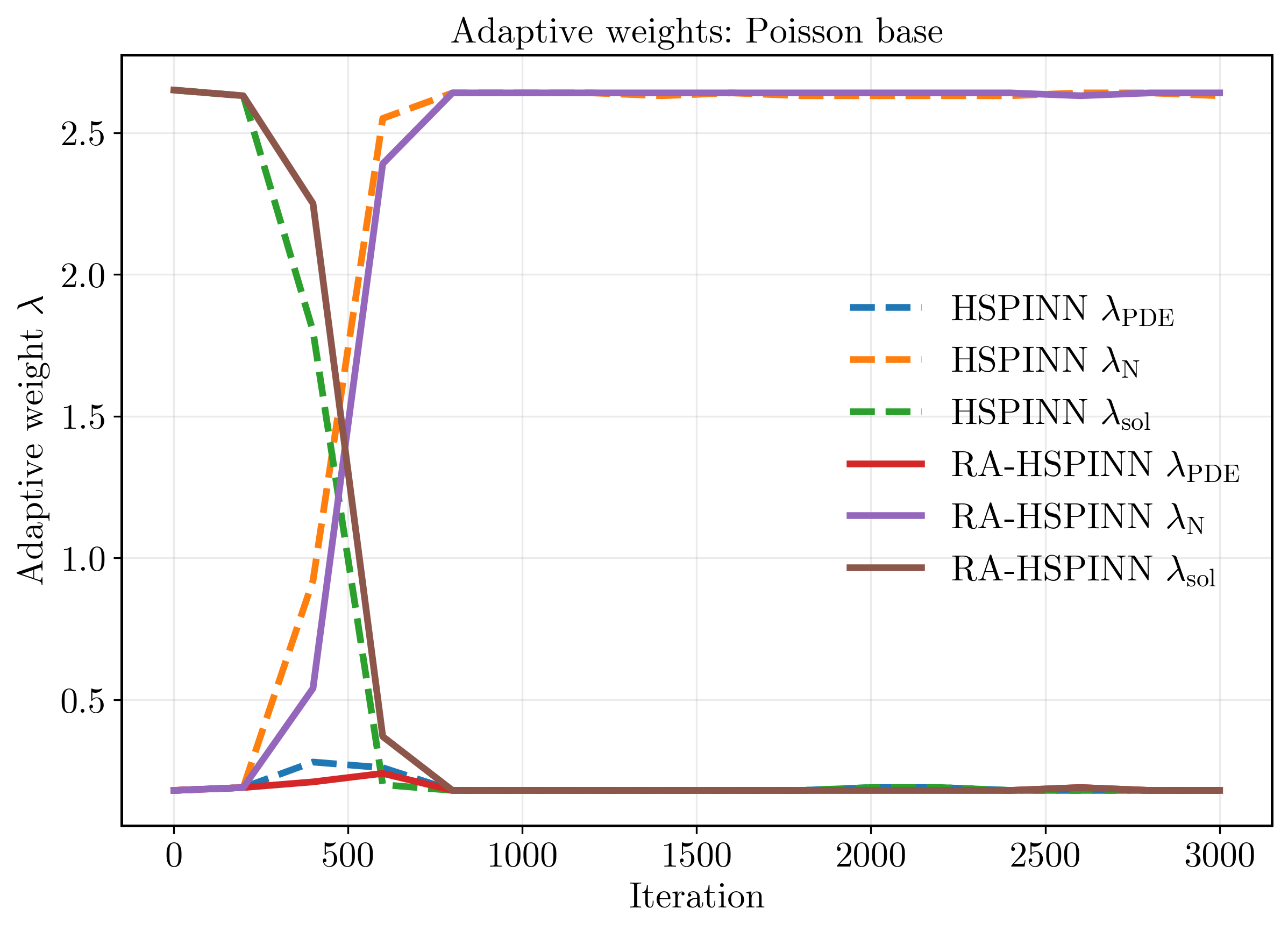}
\caption{Mixed-boundary elliptic Poisson benchmark. Adaptive global weights for
the Poisson residual, Neumann residual, and solution regularization term.}
\label{fig:fields_poisson_base_2}
\end{figure}

The adaptive-weight curves in Fig.~\ref{fig:fields_poisson_base_2} show how the
inverse-EMA rule balances the Poisson residual, the Neumann residual, and the
small solution-regularization term during training. Unlike the sharp-gradient
Burgers and mixed first-order Poisson cases, this benchmark does not exhibit a
severe residual-competition pattern: the hard mask already enforces the two
Dirichlet edges exactly, and the remaining optimization mainly adjusts the
smooth interior field while satisfying the Neumann constraints. The weight
evolution therefore supports the interpretation that the mixed-boundary Poisson
case is a moderate-improvement regime rather than a failure-recovery regime.
Together with the field and reliability-map results in
Fig.~\ref{fig:fields_poisson_base}, these curves show that RA-HSPINN improves
the smooth elliptic solution, but the gain remains modest because the baseline
hard--soft ansatz is already well aligned with the target solution.

\subsection{Mixed first-order periodic Poisson system}

The most informative elliptic-type example is the mixed first-order Poisson
system. Instead of differentiating a reliability-modulated scalar field through
a direct second-order Laplacian, the problem is written as
\begin{align}
    u_x-p&=0, \\
    u_y-q&=0, \\
    p_x+q_y&=f(x,y).
    \label{eq:mixed_first_order}
\end{align}
The neural network outputs three channels, $[u,p,q]$. In HSPINN, the output is
given by $v_\theta(\varphi(x,y))$, while in RA-HSPINN the same multi-output
field is modulated by a scalar reliability field:
\begin{equation}
    \widehat{\bm{z}}_{\theta,\phi}(x,y)
    =
    \rho_\phi(\varphi(x,y))
    \begin{bmatrix}
    u_\theta(\varphi(x,y)) \\
    p_\theta(\varphi(x,y)) \\
    q_\theta(\varphi(x,y))
    \end{bmatrix}.
    \label{eq:mixed_output}
\end{equation}
The periodic feature vector is
\begin{equation}
    \varphi(x,y)=\left[\sin x,\cos x,\sin y,\cos y\right].
\end{equation}
The optimized interior loss for this mixed system is
\begin{equation}
    \mathcal{L}_{\pde}^{\mathrm{mix}}
    =
    \mathbb{E}\!\left[(p_x+q_y-f)^2\right]
    +\mathbb{E}\!\left[(u_x-p)^2\right]
    +\mathbb{E}\!\left[(u_y-q)^2\right].
    \label{eq:mixed_pde_loss}
\end{equation}
Thus, the balance equation and the first-order consistency equations are
optimized as ordinary mean-square residuals within the same system block. This
choice avoids introducing separate manual weights for the auxiliary variables
$p$ and $q$, while the inverse-EMA rule balances the mixed PDE block against
the boundary or data loss.

The exact multi-mode solution is
\begin{equation}
    u^*(x,y)=\sin x+\sin(3x)\cos(3y)+0.5\sin(5x).
    \label{eq:mixed_multimode_exact}
\end{equation}
The forcing term is obtained from the manufactured solution by automatic
differentiation in the experiment code. With the sign convention used in
Eq.~\eqref{eq:mixed_first_order}, the system corresponds to the periodic
equation $\Delta u=f$, whereas the mixed-boundary Poisson benchmark above uses
the convention $-\Delta u=f$. The two benchmarks are therefore stated according
to their respective implementations.

This formulation is important because it clarifies the interaction between PDE
form and reliability-aware modulation. If the direct second-order operator is
applied to a reliability-modulated scalar field, the Laplacian contains the
product-rule terms
\begin{equation}
    \Delta(\rho v)=\rho\Delta v+2\nabla\rho\cdot\nabla v+v\Delta\rho.
    \label{eq:product_rule}
\end{equation}
These terms are mathematically correct, but they can increase optimization
stiffness when the reliability field varies during training. The mixed
first-order formulation avoids differentiating the product $\rho_\phi v_\theta$
twice. All residual equations in Eq.~\eqref{eq:mixed_first_order} require only
first-order derivatives of the network outputs. As a result, the reliability
field can modulate the representation without introducing a direct second-order
product-rule burden into the residual.

\begin{figure}[t]
\centering
\includegraphics[width=0.98\textwidth]{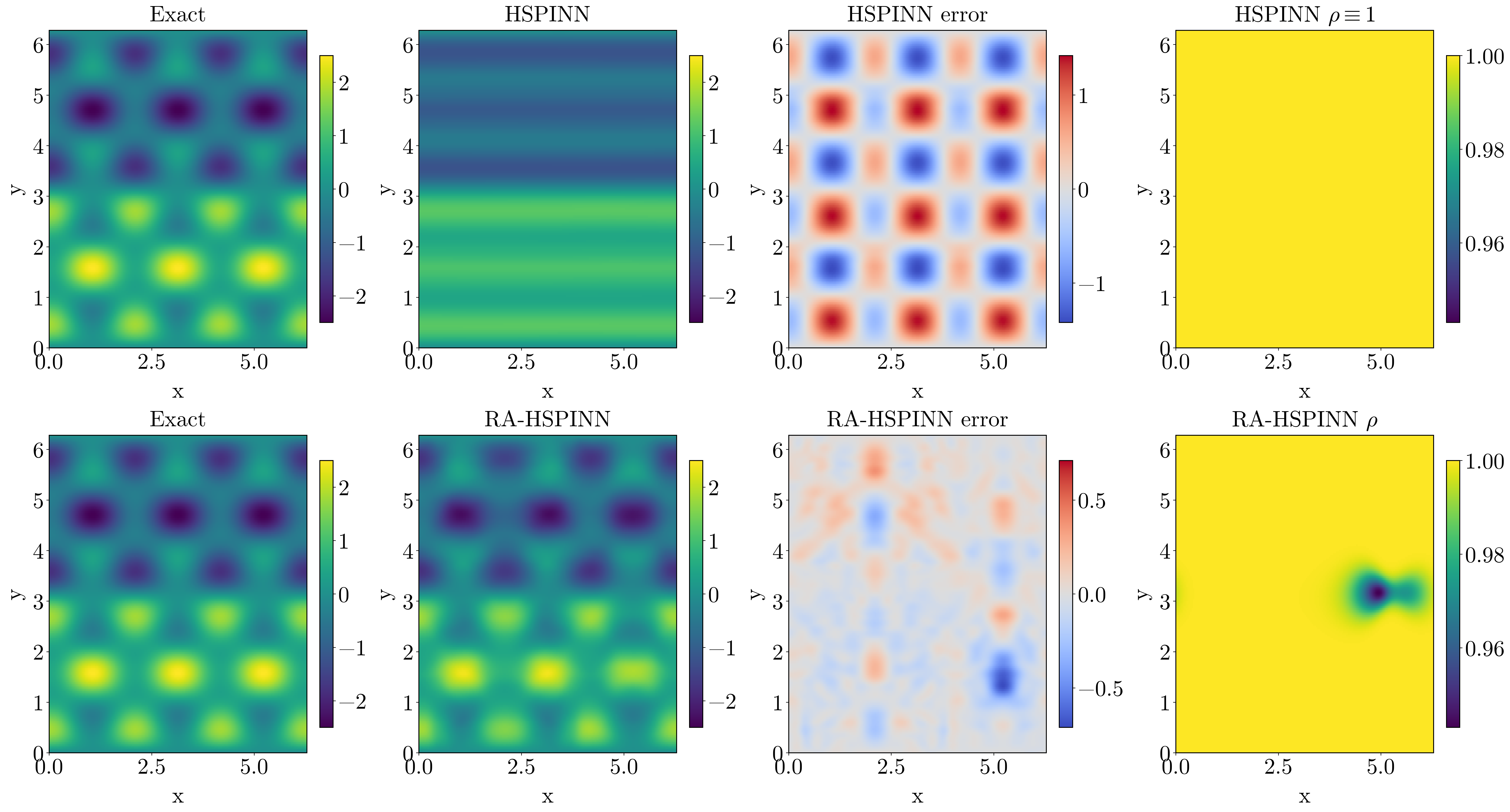}
\caption{Mixed first-order periodic Poisson benchmark with a multi-mode
solution. Reference field, HSPINN and RA-HSPINN predictions, pointwise error
fields, and the learned reliability map. RA-HSPINN captures the multi-mode
structure more accurately and substantially reduces the large error observed in
the fixed HSPINN baseline.}
\label{fig:fields_mixed_multi_mode}
\end{figure}

The numerical fields in Fig.~\ref{fig:fields_mixed_multi_mode} show a clear
difference between the two models. The fixed HSPINN baseline fails to recover
the multi-mode structure of the solution and leaves large pointwise errors over
the domain. In contrast, RA-HSPINN captures the dominant low- and
high-frequency components more accurately and substantially reduces the error
field. The learned reliability map in the same figure indicates that the
modulation is not used as a strongly suppressive mask; instead, it provides a
bounded adaptive degree of freedom that helps the multi-output representation
escape the stagnation observed in the fixed HSPINN baseline.

Quantitatively, HSPINN stagnates at a final relative error of $0.616051$,
whereas RA-HSPINN reaches $0.109856$, corresponding to an $82.17\%$ reduction.
The final RA reliability statistics are mean $0.99899$, minimum $0.94436$, and
maximum $1.0$. These values show that the reliability field remains close to
one at convergence, but the training trajectory and final solution differ
substantially from the fixed HSPINN result. Therefore, the benefit of RA-HSPINN
should not be interpreted only from the final mean reliability value; it also
depends on the optimization path created by reliability modulation and adaptive
global balancing.

\begin{figure}[t]
\centering
\includegraphics[width=0.65\textwidth]{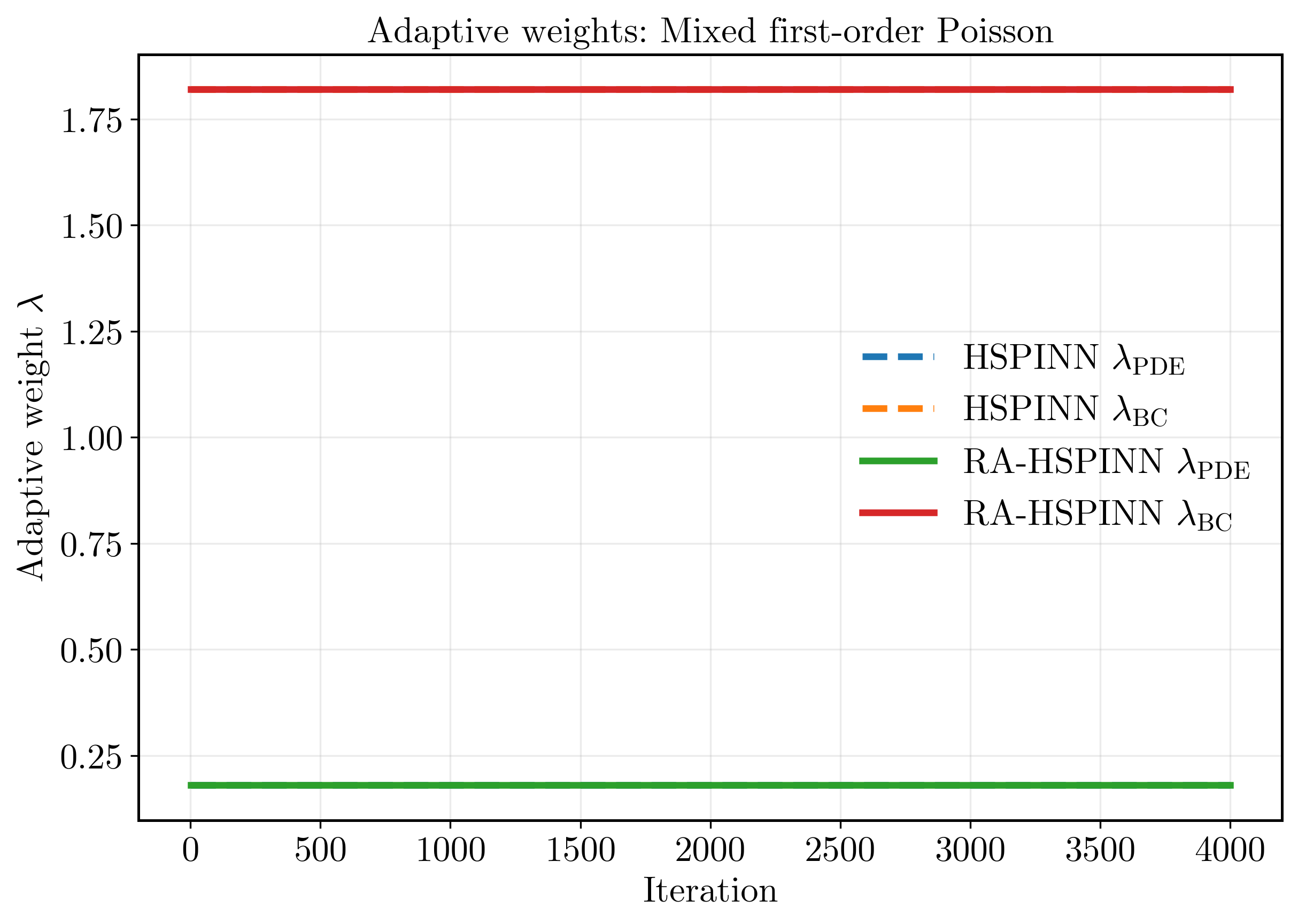}
\caption{Mixed first-order periodic Poisson benchmark with a multi-mode
solution. Adaptive global weights for the mixed PDE system block and the
boundary or data loss.}
\label{fig:fields_mixed_multi_mode_2}
\end{figure}

The adaptive-weight curves in Fig.~\ref{fig:fields_mixed_multi_mode_2} show how
the inverse-EMA rule balances the mixed PDE system block against the
boundary/data loss during training. In this benchmark, the PDE block contains
both the balance equation and the first-order consistency residuals,
\((p_x+q_y-f)\), \((u_x-p)\), and \((u_y-q)\). The weight evolution therefore
reflects the competition between satisfying the coupled first-order system and
anchoring the multi-mode solution through the boundary/data term. Together with
the field and reliability-map results in
Fig.~\ref{fig:fields_mixed_multi_mode}, these curves support the interpretation
that RA-HSPINN benefits from the combination of reliability-aware modulation,
inverse-EMA global balancing, and the first-order residual formulation. This
combination helps the optimizer escape the stagnation observed in the fixed
HSPINN baseline for the multi-mode elliptic-type system.

\subsection{Baseline choice and result coverage}
\label{sec:baseline_coverage}

The primary comparison in this paper is between RA-HSPINN and the fixed
hard--soft HSPINN baseline. This comparison directly tests the proposed
mechanism, because it evaluates whether reliability-aware modulation improves
an already admissible hard--soft neural solver. A supplementary fully soft PINN
baseline, denoted as SPINN, is also included to diagnose when hard constraint
embedding or periodic feature embedding itself is beneficial or detrimental.
Thus, SPINN is used as a diagnostic reference for the effect of hard embedding,
whereas the main interpretation of the proposed RA-HSPINN method is made
relative to HSPINN.

\subsubsection{Comparison between SPINN, HSPINN, and RA-HSPINN}
\label{sec:spinn_hspinn_comparison}

For completeness, \cref{tab:spinn_hspinn_comparison} reports a three-way
comparison between the fully soft PINN baseline (SPINN), the fixed hard--soft
PINN baseline (HSPINN), and the proposed RA-HSPINN. This comparison separates
two effects. The SPINN--HSPINN columns quantify the effect of embedding known
Dirichlet or periodic structure into the neural trial space. The HSPINN--RA
columns quantify the additional effect of reliability-aware modulation after
the hard--soft structure is already present. The reported percentages are
computed as
\[
    \frac{\varepsilon_{\mathrm{old}}-\varepsilon_{\mathrm{new}}}
    {\varepsilon_{\mathrm{old}}}\times 100\% ,
\]
so positive values indicate a reduction in relative error and negative values
indicate degradation.

\begin{table}[t]
\centering
\caption{Three-way comparison between SPINN, HSPINN, and RA-HSPINN. All models
use the same solution-network architecture (width 64, depth 4), optimizer
(Adam, $\mathrm{lr}=10^{-3}$), and collocation budget
($N_{\mathrm{int}}=2048$, 4000 steps). RA-HSPINN additionally uses the
reliability network described in \cref{sec:setup}. Positive percentages denote
relative error reduction, while negative percentages denote degradation.}
\label{tab:spinn_hspinn_comparison}
\resizebox{\textwidth}{!}{%
\begin{tabular}{lcccccc}
\toprule
Benchmark
& SPINN $L_2$
& HSPINN $L_2$
& RA-HSPINN $L_2$
& SPINN$\to$HSPINN (\%)
& HSPINN$\to$RA (\%)
& Main observation \\
\midrule
Burgers B (sharp gradient)
& $1.79\times10^{-3}$
& $3.67849\times10^{-1}$
& $4.965\times10^{-3}$
& $-20450.2$
& $98.65$
& Hard mask hurts; RA largely compensates \\

Burgers C (noisy IC)
& $7.19\times10^{-4}$
& $1.733\times10^{-3}$
& $4.78\times10^{-4}$
& $-141.03$
& $72.42$
& RA-HSPINN best under unreliable IC \\

Convection $\beta=1$
& $5.67\times10^{-3}$
& $4.432\times10^{-3}$
& $1.720\times10^{-3}$
& $21.83$
& $61.18$
& Fourier features help; RA further improves \\

Convection (steep front)
& $4.88\times10^{-1}$
& $3.309\times10^{-3}$
& $1.323\times10^{-3}$
& $99.32$
& $60.02$
& SPINN fails without periodic embedding \\

Poisson (mixed BC)
& $1.14\times10^{-1}$
& $3.5212\times10^{-2}$
& $2.4872\times10^{-2}$
& $69.11$
& $29.36$
& Hard mask dominates; RA adds modest gain \\

Mixed Poisson (multi-mode)
& $6.25\times10^{-1}$
& $6.16051\times10^{-1}$
& $1.09856\times10^{-1}$
& $1.43$
& $82.17$
& RA recovers when both baselines stagnate \\
\bottomrule
\end{tabular}}
\end{table}

\subsubsection{Analysis and insights}

The three-way comparison in \cref{tab:spinn_hspinn_comparison} shows that
neither SPINN nor HSPINN is uniformly superior. The relative ranking depends on
whether the embedded structure is well aligned with the solution and with the
resulting optimization landscape.

\paragraph{Effect of hard or periodic embedding}
The clearest positive case for embedded structure is the convection benchmark
with a steep localized front. SPINN gives a large relative error of
$4.88\times10^{-1}$, whereas HSPINN reduces the error to
$3.309\times10^{-3}$, corresponding to a $99.32\%$ reduction. This improvement
is mainly due to the periodic feature embedding, which represents spatial
periodicity directly instead of requiring the network to learn periodic
structure from raw coordinates. A similar but less dramatic trend is observed
in the mixed-boundary Poisson benchmark. There, the hard mask $xy$ enforces the
two homogeneous Dirichlet edges exactly and reduces the error from
$1.14\times10^{-1}$ to $3.5212\times10^{-2}$.

However, hard embedding is not automatically beneficial. In the sharp-gradient
Burgers benchmark, SPINN achieves $1.79\times10^{-3}$, while HSPINN stagnates at
$3.67849\times10^{-1}$. This result indicates that the boundary-vanishing hard
mask can make the fixed hard--soft representation poorly conditioned for this
high-frequency manufactured solution. The noisy Burgers benchmark shows a
milder version of the same issue: HSPINN is worse than SPINN because the model
must reconcile the hard-constrained representation with corrupted and
boundary-incompatible initial data. These cases show why SPINN is useful as a
diagnostic baseline even though the central comparison of the paper is between
HSPINN and RA-HSPINN.

\paragraph{Effect of reliability-aware modulation}
The HSPINN-to-RA columns isolate the contribution of the proposed
reliability-aware extension. RA-HSPINN improves over HSPINN in all six reported
benchmarks. The relative error reductions range from $29.36\%$ in the smooth
mixed-boundary Poisson problem to $98.65\%$ in the sharp-gradient Burgers
problem. The largest gains occur when the fixed hard--soft representation
stagnates or becomes poorly conditioned. In Burgers B, RA-HSPINN reduces the
error from $3.67849\times10^{-1}$ to $4.965\times10^{-3}$. In the mixed
first-order Poisson problem, it reduces the error from $6.16051\times10^{-1}$
to $1.09856\times10^{-1}$.

The results also show that the reliability field is not merely a small
regularization add-on. In the difficult Burgers and mixed Poisson cases, the
fixed HSPINN baseline fails to reach an accurate solution, while the
reliability-aware representation recovers a much lower error without relaxing
the embedded structure. This supports the main hypothesis of the paper:
RA-HSPINN is most useful when the hard--soft trial space is admissible but
difficult to optimize with a fixed interior representation.

\paragraph{Scope of the method}
The comparison also clarifies the intended scope of RA-HSPINN. The method is
not claimed to replace SPINN in every setting. In the sharp-gradient Burgers
benchmark, SPINN still gives a smaller error than RA-HSPINN because the
hard-constrained representation is particularly unfavorable for that case.
Instead, the role of RA-HSPINN is to improve a hard--soft solver while
preserving its embedded constraints. For this reason, the primary comparison in
this paper remains HSPINN versus RA-HSPINN, while the SPINN results are
included mainly to reveal when hard embedding itself is helpful or harmful.

\paragraph{Mixed first-order Poisson}
The mixed first-order Poisson problem is especially informative because both
SPINN and HSPINN struggle with the multi-mode solution, giving relative errors
of $6.25\times10^{-1}$ and $6.16051\times10^{-1}$, respectively. RA-HSPINN
reduces the error to $1.09856\times10^{-1}$, corresponding to an $82.17\%$
improvement over HSPINN. This result supports the interpretation developed in
the preceding subsection: reliability-aware modulation is especially effective
when combined with a residual formulation that allows the optimizer to escape
the stagnation of the fixed baseline representation.

\subsection{Ablation studies}
\label{sec:ablation}

The preceding comparison separates the fully soft PINN baseline, the fixed
hard--soft HSPINN baseline, and the final RA-HSPINN model. To further clarify
which component of the proposed method is responsible for the observed gains,
\cref{tab:ablation_components} reports a controlled component ablation over the
mechanisms retained in the final formulation. Exploratory local-reweighting
variants are intentionally omitted because they are outside the final
RA-HSPINN formulation. The ablation therefore focuses on two retained
ingredients: inverse-EMA global loss balancing and the learned reliability
field with its regularization.

The absolute errors in this ablation should be interpreted within the controlled
ablation setting, because the purpose is to compare model components under the
same training budget, collocation budget, architecture, and seed. Thus, the
most important comparison is the relative behavior among the three rows in
\cref{tab:ablation_components}, rather than exact numerical equality with the
main three-way comparison in \cref{tab:spinn_hspinn_comparison}.

\begin{table}[t]
\centering
\caption{Component ablation for the mechanisms retained in the final
RA-HSPINN formulation. All variants use the same solution-network size,
optimizer, training budget, collocation budget, and seed. The second row adds
inverse-EMA global loss balancing to HSPINN. The final row adds the bounded
reliability field and its regularization on top of the global weighting
mechanism.}
\label{tab:ablation_components}
\resizebox{\textwidth}{!}{%
\begin{tabular}{lcccccccc}
\toprule
Variant
& Reliability field
& Global weights
& $\rho$ reg.
& Burgers B
& Burgers C
& Conv. steep
& Poisson
& Mixed Poisson \\
\midrule
HSPINN
& no
& no
& no
& $2.64\times10^{-1}$
& $1.52\times10^{-3}$
& $6.49\times10^{-3}$
& $1.16\times10^{-1}$
& $5.39\times10^{-1}$ \\

HSPINN + global weights
& no
& yes
& no
& $3.80\times10^{-1}$
& $5.26\times10^{-4}$
& $7.31\times10^{-3}$
& $5.21\times10^{-2}$
& $6.16\times10^{-1}$ \\

RA-HSPINN
& yes
& yes
& yes
& $\mathbf{6.64\times10^{-3}}$
& $\mathbf{2.76\times10^{-4}}$
& $\mathbf{4.32\times10^{-3}}$
& $\mathbf{2.42\times10^{-2}}$
& $\mathbf{6.82\times10^{-2}}$ \\
\bottomrule
\end{tabular}}
\end{table}

The ablation results provide a more precise interpretation of the proposed
method than the full-model comparison alone. First, the reliability field is
the dominant component in the most difficult cases. In the sharp-gradient
Burgers benchmark, adding inverse-EMA global weights alone does not help; the
error increases from $2.64\times10^{-1}$ to $3.80\times10^{-1}$. In contrast,
the full RA-HSPINN reduces the error to $6.64\times10^{-3}$, corresponding to a
$97.48\%$ reduction relative to the bare HSPINN baseline. This confirms that
global balancing alone is not sufficient when the fixed hard--soft
representation itself is poorly conditioned.

Second, inverse-EMA global loss balancing is useful in some regimes but is not
uniformly beneficial. It improves the noisy Burgers case from
$1.52\times10^{-3}$ to $5.26\times10^{-4}$ and the mixed-boundary Poisson case
from $1.16\times10^{-1}$ to $5.21\times10^{-2}$. However, it worsens the
sharp-gradient Burgers case, the steep-front convection case, and the mixed
first-order Poisson case. This shows that adaptive global weighting should be
viewed as a stabilizing mechanism for heterogeneous loss terms, not as a
complete solution to all approximation and conditioning difficulties.

Third, the full RA-HSPINN is the best retained variant in all five ablation
benchmarks. It reduces the relative error by approximately $97.48\%$ for the
sharp-gradient Burgers problem, $81.84\%$ for the noisy Burgers problem,
$33.44\%$ for the steep-front convection problem, $79.14\%$ for the
mixed-boundary Poisson problem, and $87.35\%$ for the mixed first-order Poisson
problem, all relative to the bare HSPINN row. The improvement is largest when
the fixed hard--soft representation stagnates, as in the sharp-gradient Burgers
and mixed first-order Poisson examples. The gain is more moderate for the
steep-front convection problem, where the periodic feature embedding already
captures much of the transport structure.

Finally, the ablation supports the intended positioning of RA-HSPINN. The
method should not be interpreted as merely an adaptive loss-weighting strategy.
Global weighting can help when the main difficulty is imbalance among residual
terms, but the reliability field provides an additional representational
mechanism that is decisive when global weights alone are insufficient. This
distinction is important for the proposed contribution: RA-HSPINN combines
loss-level balancing with a reliability-modulated hard--soft trial space, and
the latter is responsible for the strongest recoveries observed in the
controlled ablation.

\section{Discussion}
\label{sec:discussion}

The numerical results support a targeted interpretation of RA-HSPINN. The
method is most effective when the fixed hard--soft trial space is admissible
but difficult to optimize with a single unmodulated free component. The
sharp-gradient Burgers benchmark and the mixed first-order multi-mode Poisson
benchmark provide the clearest evidence. In the sharp-gradient Burgers case,
HSPINN stagnates at a relative error of $3.67849\times10^{-1}$, whereas
RA-HSPINN reaches $4.965\times10^{-3}$, corresponding to a $98.65\%$ reduction.
In the mixed first-order Poisson case, RA-HSPINN reduces the error from
$6.16051\times10^{-1}$ to $1.09856\times10^{-1}$, an $82.17\%$ improvement.
These results show that reliability-aware modulation can substantially improve
a hard--soft representation without relaxing the embedded boundary or periodic
structure.

The Burgers benchmark with noisy and incompatible initial data shows a
different but related benefit. Both models use the same hard Dirichlet ansatz
and the same hand-designed weighting of the unreliable initial data. Therefore,
the improvement from $1.733\times10^{-3}$ for HSPINN to $4.78\times10^{-4}$ for
RA-HSPINN is not caused by a different boundary treatment or a different manual
data weighting. Instead, the learned reliability field and adaptive global loss
weights help the model reconcile the clean PDE residual with corrupted initial
information. The final adaptive weights shift toward the PDE residual, which is
consistent with the benchmark design because the initial condition
intentionally contains noise and boundary-incompatible perturbations.

The convection and mixed-boundary Poisson results give a more nuanced picture.
For smooth periodic convection, the final mean reliability remains close to one,
indicating that the reliability field is nearly inactive when the periodic
feature embedding already provides a well-conditioned representation. The
localized-front convection case also has reliability values close to one on
average, but RA-HSPINN still improves the relative error from
$3.309\times10^{-3}$ to $1.323\times10^{-3}$. This suggests that mild
modulation can be sufficient when the dominant periodic structure is already
encoded but the transported solution contains a localized peak. In the
mixed-boundary Poisson benchmark, the gain is more moderate ($29.36\%$ relative
to HSPINN), which is expected because the solution is smooth and the hard mask
$xy$ already enforces the dominant Dirichlet constraints exactly.

The three-way comparison with SPINN clarifies that hard embedding is problem
dependent. In the steep-front convection and mixed-boundary Poisson benchmarks,
hard or periodic embedding substantially improves over the fully soft baseline.
However, in the sharp-gradient Burgers benchmark, the fully soft SPINN baseline
achieves a lower error than RA-HSPINN. This does not contradict the proposed
method, because RA-HSPINN is not intended to replace SPINN in every setting.
Its role is to improve an already hard--soft solver while preserving its
embedded constraints. The SPINN baseline is therefore best interpreted as a
diagnostic tool for understanding when hard embedding itself is helpful or
harmful.

The ablation further separates the two retained mechanisms:
inverse-EMA global loss balancing and reliability-aware modulation. Global
balancing alone improves some cases, such as noisy Burgers and mixed-boundary
Poisson, but it is not uniformly beneficial and can worsen other cases. In
contrast, the full RA-HSPINN is the best retained variant in all five ablation
benchmarks. Relative to bare HSPINN, it reduces the error by approximately
$97.48\%$ for sharp-gradient Burgers, $81.84\%$ for noisy Burgers, $33.44\%$
for steep-front convection, $79.14\%$ for mixed-boundary Poisson, and
$87.35\%$ for mixed first-order Poisson. These ablation results indicate that
the reliability field is not merely a regularization add-on; it provides an
additional representational mechanism that becomes decisive when loss-level
balancing alone is insufficient.

The mixed first-order Poisson result also suggests an important formulation
principle. Reliability-aware modulation interacts with the differential order
of the PDE operator. If a direct second-order operator is applied to a product
$\rho_\phi v_\theta$, product-rule terms involving derivatives of $\rho_\phi$
appear in the residual. These terms are mathematically correct but can increase
optimization stiffness. The mixed first-order formulation avoids differentiating
the reliability-modulated product twice, so the reliability field can modulate
the representation without introducing a direct second-order product-rule
burden. This observation is important for computational physics applications:
the success of a physics-informed neural method depends not only on the network
architecture, but also on how the continuous operator is formulated before
training.

It is worth noting that the Burgers residual also involves a second-order derivative applied to the RA-HSPINN ansatz $x(1-x)\rho_\phi v_\theta$. Expanding $\partial_{xx}[x(1-x)\rho_\phi v_\theta]$ yields product-rule terms involving $\partial_x\rho_\phi$, $\partial_{xx}\rho_\phi$, and the second derivative of the mask $x(1-x)$. In the one-dimensional Burgers case these terms are bounded because the mask is fixed and the reliability smoothness regularization (Eq.~\eqref{eq:rho_reg}, $\alpha_s = 10^{-5}$) discourages large spatial derivatives of $\rho_\phi$. The product-rule burden is therefore present but controlled. In contrast, for a two-dimensional Laplacian acting
on $\rho_\phi v_\theta$ with no mask and no smoothness constraint, the interaction between $\nabla\rho_\phi$ and $\nabla v_\theta$ in both spatial directions creates a larger and less controlled optimization stiffness. The mixed first-order reformulation eliminates this entirely by ensuring that no second derivative of the modulated product enters the residual.

The computational cost of RA-HSPINN is higher than that of HSPINN because the
method introduces an auxiliary reliability network, additional forward
evaluations, and derivative terms for reliability regularization. This overhead
should be considered when the base HSPINN already solves a smooth problem to a
satisfactory tolerance. The practical value of RA-HSPINN is therefore
problem-dependent. It is most justified when the baseline hard--soft solver
stagnates, when the residual structure is localized or multi-modal, or when the
training constraints contain unreliable components that should not dominate the
PDE residual.

\section{Conclusion}
\label{sec:conclusion}

This paper introduced RA-HSPINN, a reliability-aware extension of hard--soft
physics-informed neural networks. The method preserves embedded Dirichlet or
periodic constraints while adding a bounded learnable reliability field to the
admissible interior representation. It also uses inverse-EMA global loss
balancing to reduce domination among heterogeneous residual terms. The
reliability field is not interpreted as a physical parameter or calibrated
uncertainty; it is a numerical modulation variable that enriches the hard--soft
trial space without relaxing the prescribed constraints.

The method was evaluated on nonlinear Burgers equations, periodic convection, a
mixed-boundary Poisson problem, and a mixed first-order Poisson system. Across
the main HSPINN-to-RA comparison, RA-HSPINN improved all six reported
benchmarks, with reductions of $98.65\%$ for sharp-gradient Burgers,
$72.42\%$ for noisy and incompatible Burgers initial data, $61.19\%$ for smooth
periodic convection, $60.02\%$ for localized-front convection, $29.36\%$ for
mixed-boundary Poisson, and $82.17\%$ for mixed first-order Poisson. The
controlled ablation confirms that inverse-EMA global weighting is helpful in
some regimes but is not sufficient by itself; the reliability field provides
the strongest recovery when the fixed hard--soft representation stagnates.

The main numerical conclusion is conditional rather than universal. RA-HSPINN
is most useful when a hard--soft trial space is admissible but poorly
conditioned, especially for localized, unreliable-data, or multi-mode residual
structures. It should be viewed as a targeted extension of HSPINN, not as a
replacement for all soft or hard-constrained PINN formulations. The mixed
first-order Poisson result further suggests that reliability-aware neural
modulation should be designed together with the PDE formulation, since
first-order residual structures can avoid unnecessary product-rule stiffness.
Future work should test the method on higher-dimensional PDEs, complex
geometries, repeated random seeds, and adaptive sampling strategies coupled to
the learned reliability field.

\section*{Acknowledgments}
This research was funded by the National Foundation for Science and Technology Development (NAFOSTED) through project 107.02-2025.21. The authors are grateful for the support.

\section*{CRediT authorship contribution statement}
\textbf{Duc Tien Nguyen:} Writing -- original draft, Visualization, Validation, Software, Methodology, Investigation, Formal analysis, Data curation, Conceptualization. \textbf{Hang Tran:} Writing -- original draft, Visualization, Validation, Software, Methodology, Investigation, Formal analysis, Data curation.
\textbf{Trinh Minh Tuan:} Writing -- original draft, Data curation, Validation, Investigation.
\textbf{Nguyen Duc Manh:} Writing -- review \& editing, Formal analysis.
\textbf{Dinh Gia Ninh:} Writing -- review \& editing, Methodology, Conceptualization, Investigation, Supervision.

\section*{Declaration of competing interest}
The authors declare that they have no known competing financial interests or personal relationships that could have appeared to influence the work reported in this paper.


\bibliography{Bio}

\appendix

\end{document}